\documentclass[11pt]{article}
\usepackage{amsmath,amssymb,amsfonts,amsthm}
\usepackage{colordvi}
\usepackage[pdftex]{graphicx}
\usepackage[pdftex]{color}
\usepackage{epstopdf}
\usepackage{verbatim}    
\usepackage{array}       
\usepackage{algorithmic}
\usepackage{algorithm}
\usepackage{enumerate}
\usepackage{mdframed}

\let\counterwithin\relax
\usepackage{chngcntr}
\counterwithin{equation}{section}

\usepackage[sort]{cite}

\mdfsetup{
  linewidth=.9bp,
  innerleftmargin=5bp,
  innerrightmargin=5bp,
  innertopmargin=5bp,
  innerbottommargin=5bp,
}

\def\bx{\mathbf{x}}
\def\by{\mathbf{y}}
\def\rhoydy{\rho(\by)d\by}
\def\bu{\mathbf{u}}
\def\bv{\mathbf{v}}

\def\bm{\mathbf{m}}
\def\bn{\mathbf{n}}
\def\bp{\mathbf{p}}

\newcommand{\btau} { \mbox{\boldmath $\tau$} }
\def\R{\mathbb{R}}
\def\grad{\nabla}
\def\Gm{\Gamma}
\def\Om{\Omega}
\def\div{\nabla\cdot}

\def\<{\langle}
\def\>{\rangle}
\def\S{\mathbb{S}}

\def\nomega{N_\Omega}
\def\nterm{N_{\textrm{term}}}  
\def\ncoll{N_{\textrm{coll}}}
\def\nreal{N_{\textrm{real}}}

\def\ndofi{N_{\textrm{dof},i}}

\def\cgiter{N_{\textrm{iter}}}
\def\lmax{\ell_{\textrm{max}}}
\def\lmin{\ell_{\textrm{min}}}

\def\k{\{k\}}

\def\bl{{\boldsymbol\lambda}}
\def\d{\partial}
\def\ff{\mathbf{f}}
\def\Rset{\mathbb{R}}
\def\bT{{\bf T}}
\def\bD{{\bf D}}
\def\n{{\bf n}}

\begin{document}

\title{Stochastic multiscale flux basis for Stokes-Darcy 
flows}


\author{Ilona Ambartsumyan\thanks{Department of Mathematics, University of
Pittsburgh, 301 Thackeray Hall, Pittsburgh, Pennsylvania 15260, USA;~{\tt \{ila6@pitt.edu, elk58@pitt.edu, chw92@pitt.edu, yotov@math.pitt.edu\}}; partially supported by NSF grants DMS 1418947 and DMS 1818775 and DOE grant DE-FG02-04ER25618.}~\and
Eldar Khattatov\footnotemark[1]\and
ChangQing Wang\footnotemark[1]\and
Ivan Yotov\footnotemark[1]}

\date{\today}
\maketitle

\begin{abstract}
Three algorithms are developed for uncertainty quantification in
modeling coupled Stokes and Darcy flows. The porous media may consist
of multiple regions with different properties. The permeability is
modeled as a non-stationary stochastic variable, with its log
represented as a sum of local Karhunen-Lo\`eve (KL) expansions.  The
problem is approximated by stochastic collocation on either
tensor-product or sparse grids, coupled with a multiscale mortar mixed
finite element method for the spatial discretization. A
non-overlapping domain decomposition algorithm reduces the global
problem to a coarse scale mortar interface problem, which is solved by
an iterative solver, for each stochastic realization. In the
traditional domain decomposition implementation, each subdomain solves
a local Dirichlet or Neumann problem in every interface iteration. To
reduce this cost, two additional algorithms based on deterministic or
stochastic multiscale flux basis are introduced. The basis consists of
the local flux (or velocity trace) responses from each mortar degree
of freedom. It is computed for each subdomain independently before the
interface iteration begins. The use of the multiscale flux basis
avoids the need for subdomain solves on each iteration.  The
deterministic basis is computed at each stochastic collocation and
used only at this realization. The stochastic basis is formed by
further looping over all local realizations of a subdomain's KL region
before the stochastic collocation begins. It is reused over multiple
realizations. Numerical tests are presented to illustrate the
performance of the three algorithms, with the stochastic multiscale
flux basis showing significant savings in computational cost compared 
to the other two algorithms. 

\end{abstract}

\begin{keywords}
 stochastic collocation, domain decomposition, multiscale basis,
 mortar finite element, mixed finite element, Stokes-Darcy flows
\end{keywords}


\section{Introduction} \label{intro}

Coupled Stokes-Darcy flows arise in numerous applications, including
interaction between surface and groundwater flows, cardiovascular
flows, industrial filtration, fuel cells, and flows in fractured or
vuggy reservoirs. The Stokes equations describe the motion of
incompressible fluids and the Darcy model describes the infiltration
process. In this work we consider mixed velocity-pressure formulations
in both regions. The Beavers-Joseph-Saffman slip with friction
condition \cite{BeaverJoseph,Saffman} is applied on the Stokes-Darcy
interface.  Existence and uniqueness of a weak solution has been
proved in \cite{LaytonSchieweckYotov}, see also \cite{DMQ} for
analysis of the system with a pressure Darcy formulation. The finite
element approximation of the coupled problem has been studied
extensively, see, e.g.,
\cite{LaytonSchieweckYotov,DMQ,RivYot,Galvis-Sarkis} for some of the
early works. To the best of our knowledge, all previous studies have
considered the deterministic model. Often, due to incomplete knowledge
of the physical parameters, quantification of the model uncertainty
needs to be incorporated. In this paper, we study the interaction of a
free fluid with a porous medium with stochastic permeability. Even
though the stochasticity comes only from the uncertain nature of the
permeability in the porous region, the resulting solution is
stochastic over the entire domain, due to the coupling conditions. We
develop algorithms for computing statistical moments of the
solution, such as mean, variance, and higher moments.

In this work the permeability function is parameterized using a truncated
Karhunen-Lo\'eve (KL) expansion with independent identically
distributed random variables as coefficients. Given a covariance
relationship with empirically determined statistics, one can compute
the eigenvalues and corresponding eigenfunctions that form the KL
series. This approach is commonly used for stochastic permeability as
shown in \cite{lu2008ssf,winter-tart-2002,Zhang_book,zhang2004eho} and
in particular can be used in the framework for stochastic collocation
and mixed finite elements \cite{ganis2008sc}. Following
\cite{lu2008ssf,GanisYotovZhong}, we consider non-stationary porous
media with different covariance functions in different parts of the
domain, which allows for modeling heterogeneous media. For instance, the
arrangement of sedimentary rocks in distinct layers motivates the use
of such statistically independent regions, each region corresponding
to a particular rock type. In this paper such regions are referred to
as KL regions. The covariance between two points that lie in different
KL regions is zero, while otherwise it depends on the distance between
those points.

To compute the statistical moments of the solution, we employ the
stochastic collocation method
\cite{babuska2008stochastic,Xiu-Hest,Nob-Temp-Web,ganis2008sc}, which
samples the stochastic space at specifically chosen points.  In this
work we use zeros of orthogonal polynomials with either tensor product
or sparse grid approximations. The method is non-intrusive, since it
requires solving a sequence of deterministic problems. In this regard
it resembles Monte Carlo Simulation (MCS) \cite{fishman1996monte}.
However, MCS may exhibit a high computational cost due to the need to
generate valid representative statistics from a large number of
realizations at random points in the stochastic event space. In
comparison, stochastic collocation provides improved approximation of
polynomial interpolation type and thus may result in better accuracy
than MCS with fewer realizations. In terms of its approximation
properties, stochastic collocation is comparable to polynomial chaos
expansions \cite{XiuandKarniadakis2002} or stochastic finite element
methods \cite{deb2001ssp,ghanem2003sf}. However, the intrusive
character of these methods may complicate their implementation. In
addition, they result in high dimensional algebraic problems with
fully coupled physical and stochastic dimensions.

In each stochastic realization we solve the coupled Stokes-Darcy
problem using the multiscale mortar mixed finite element method
(MMMFEM) introduced in \cite{arbogast2007mmm,GVY,VasWangYot}. The
MMMFEM decomposes the physical domain into a union of non-overlapping
subdomains of Stokes or Darcy type. Any Darcy subdomain is assumed to
be contained in only one KL region, but KL regions may contain
multiple subdomains. Each subdomain is discretized on a fine scale
using stable and conforming mixed finite element spaces of Stokes or
Darcy type. The grids are allowed to be non-matching along subdomain
interfaces. A mortar space is introduced on the interfaces and
discretized on a coarse scale. A coarse scale mortar Lagrange
multiplier $\lambda_H$ is used to impose weakly continuity of flux.
Since we allow for multiple subdomains, we must account for interfaces
of Stokes-Darcy, Darcy-Darcy and Stokes-Stokes types.  In particular,
on Stokes-Darcy or Darcy-Darcy interfaces $\lambda_H$ is the normal
stress or pressure, respectively - a scalar quantity, and it is used
to impose continuity of the normal velocity. On Stokes-Stokes
interfaces $\lambda_H$ is the normal stress vector, and it is used to
impose continuity of the entire velocity vector.  Following
\cite{glowinski1988dda,VasWangYot}, the global fine scale problem is
reduced to a coarse scale interface problem, which can be solved in
parallel by Krylov iterative solvers. In the traditional
implementation, the action of the interface operator requires solving
Neumann problems in Stokes subdomains and Dirichlet problems in Darcy
subdomains. The finite element tearing and interconnecting (FETI)
method \cite{FR91,TW} is employed to deal with the possibly singular
Stokes subdomain problems.  For previous work on domain decomposition
for Stokes-Darcy flows in the two-subdomain case, we refer the reader
to
\cite{Disc-Quart-2003,Disc-Quart-2004,Hoppe-Porta-Vass,Galvis-Sarkis-DD}.

The multiscale approximation of the MMMFEM is motivated
by the fact that in porous media problems, resolving fine scale
accuracy is often computationally infeasible. The method is an
alternative to other multiscale methods, such as the variational
multiscale method \cite{hughes1998variational,Arbogast_2004} and
multiscale finite elements
\cite{hou1997multiscale,Chen_Hou_2003,AKL_2006,Chung-Eff-Lee,Eff-Galvis-Hou,Galvis-Eff,Chung-Ef-Hou}.
Both have been applied to stochastic problems in
\cite{velamur2006stochastic,GZ07DataDriven} and
\cite{dostert2008multiscale,Aarn-Efen-SISC-2008} respectively. The
MMMFEM is more flexible than the aforementioned methods, since it
allows for a posteriori error estimation and adaptive refinement of
the coarse scale mortar interface mesh \cite{arbogast2007mmm}. 

In the traditional implementation of the MMMFEM, the dominant
computational cost is the solution of Stokes or Darcy subdomain
problems at each interface iteration. Even though the dimension of the
interface problem is reduced due to the coarse scale mortar space, the
cost can be significant and the number of iterations grow with the
condition number of the interface operator. To alleviate this cost, we
utilize a multiscale flux basis. We follow the approach introduced in
\cite{ganis2009ms} for the Darcy problem, extended to the Stokes-Darcy
problem in \cite{GanVasWangYot}, and the stochastic Darcy problem in
\cite{GanisYotovZhong}. In particular, we compute a multiscale basis
consisting of the local fine scale flux (or velocity trace) responses
from each coarse scale mortar degree of freedom. It is computed for
each subdomain independently before the interface iteration
begins. The use of the multiscale flux basis avoids the need for
subdomain solves on each Krylov iteration, since the action of the
interface operator can be computed by a simple linear combination of
the basis functions. We develop and compare two multiscale flux basis
algorithms for the stochastic Stokes-Darcy problem. In the first
method, referred to as deterministic multiscale flux basis, the basis
is computed at each stochastic collocation and used only during the
interface iteration at this realization. In this approach the number
of subdomain solves is proportional to the number of mortar degrees of
freedom per subdomain and the total number of stochastic collocation
points. In the second method we follow the approach from
\cite{GanisYotovZhong} and explore the fact that the stochastic
parameter is represented as a sum of local KL expansions.  Unlike the
first method, we do not recompute the basis at the beginning of every
stochastic collocation. Instead, a full multiscale basis
pre-computation is carried out on each subdomain for all local
stochastic realizations before the global collocation loop
begins. This stochastic basis can be reused multiple times during the
collocation loop, since the same local subdomain stochastic structure
occurs over multiple realizations. In this method, referred to as
stochastic multiscale flux basis, the number of
subdomain solves is proportional to the number of mortar degrees of
freedom per subdomain and the number of local stochastic collocation
points, resulting in additional significant computational savings.

While in this work we focus on single phase flow, extensions to
multiphase flow can also be considered. For example, multiscale mixed
finite element methods for multiphase flow in porous media have been
developed in \cite{Arbogast-2phase,AKL_2006,frozenJac}. Stochastic
collocation with mixed finite elements for multiphase Darcy flow have
been studied in \cite{ganis2008sc}. While multiphase Stokes-Darcy
couplings have been less extensively studied, we refer the reader to
\cite{Dawson-sw-Richards} for the coupling of the shallow water
equations with the Richards equation for modeling unsaturated
groundwater flows, as well as \cite{Sun-etal-2phaseSD,Chen-Han-Wang} for using
diffuse interface models for two-phase Stokes-Darcy flows.
  
The remainder of the paper is organized as follows. The model problem
is introduced in Section \ref{model}. The MMMFEM is introduced in
Section \ref{MMMFEM}. Tensor product and sparse grid stochastic
collocation algorithms are presented in Section
\ref{collocation}. Three algorithms, including traditional MMMFEM,
deterministic and stochastic multiscale flux basis are presented in
Section \ref{algorithms}. In Section \ref{numerical}, the algorithms
are employed in several numerical tests and compared in terms of
computational efficiency.


\section{Model problem}\label{model}

Let $\Om$ be a stochastic space with probability measure $P$. For any
random variable $\xi(\omega): \Om\rightarrow\Rset$, $\xi(w) = z$, 
with probability
density function (PDF) $\rho(z)$, its mean or expectation is defined
by
\begin{equation}
E[\xi]=\int_\Om \xi(\omega)dP(\omega)=\int_\R z\rho(z) dz,\label{exp}
\end{equation}
and its variance is given by
$$Var[\xi]=E[\xi^2]-(E[\xi])^2.$$    

We denote the fluid region and the porous media region by
$D_s\subset\Rset^d$ and $D_d\subset\Rset^d$, respectively, where
$d=2,3$. Each region may consists of several non-connected components.
Let $\Gm_s$ be the outside boundary of $D_s$ and let $\n_s$ be the
outward unit normal vector on $\partial D_s$.  Let $\Gm_d$ be the
outside boundary of $D_d$ and let $\n_d$ be the outward unit normal
vector on $\partial D_d$.  The entire physical domain is defined as
$D=D_s\cup D_d$, with the Stokes-Darcy interface denoted by
$\Gm_{sd}=D_s\cap D_d$. Let $(\bu_s, p_s)$ and $(\bu_d, p_d)$ be the
velocity and pressure unknowns in the Stokes and Darcy regions,
respectively.  In the Stokes region, let $\nu_s$ be the viscosity
coefficient and define the deformation rate tensor $\bD$ and stress
tensor $\bT$ by
$$
\mathbf {D}(\mathbf{u}_s):= \frac{1}{2}(\nabla\mathbf{u}_s +
(\nabla\mathbf{u}_s)^T ), \quad \mathbf{T}(\mathbf{u}_s,p_s): = -
p_{s}\mathbf{I} + 2\nu_s \mathbf{ D}(\mathbf{u}_s). 
$$ 

In the Darcy region, let $\nu_d$ be the viscosity coefficient and
$K(\bx,\omega)$ be a stochastic function defined on $D \times \Om$
representing the non-stationary permeability of the porous medium. We
assume $K$ to be uniformly positive definite for $P$-almost every
$\omega\in\Om$ with components in $L^{\infty}(D_d)$. The coupled
Stokes-Darcy flow model is as follows. For $P$-almost every $\omega\in\Om$, in
the Stokes region, $(\bu_s,p_s)$ satisfy
\begin{align}
-\div \bT(\bu_s,p_s) \equiv 
-2\nu_s \div \bD(\bu_s) + \grad p_s & =  \ff_s \quad \text{in } D_s, \label{stokes.1}\\
\div \bu_s & =  0  \quad \text{in } D_s, \label{stokes.2}\\
\bu_s & =  { \bf 0 } \quad \text{on } \Gm_s, \label{stokes.3}
\end{align}
where $\ff_s$ represents the body force.  In the Darcy region, $(\bu_d, p_d)$ satisfy
\begin{align}
\nu_d K(\bx,\omega)^{-1}\bu_d + \grad p_d & =  \ff_d \quad \text{in } D_d,\label{darcy.1} \\
\div \bu_d & =  q_d  \quad \text{in } D_d, \label{darcy.2}\\
\bu_d\cdot\bn_d & =  0 \quad \text{on } \Gm_d, \label{darcy.3}
\end{align}
where $\ff_d$ is the gravity force and $q_d$ is an external source or sink term satisfying the solvability condition $$\int_{D_d}q_d \  d\bx=0.$$
The two models are coupled on $\Gm_{sd}$ through the following interface conditions: 
\begin{align}
 \bu_s\cdot\bn_s + \bu_d\cdot\bn_d & =  0 \quad \text{on } \Gm_{sd}, 
 \label{cont.flux}\\
- (\bT(\bu_s,p_s) \bn_s) \cdot \bn_s \equiv
p_s - 2\nu_s (\bD(\bu_s) \bn_s) \cdot \bn_s & =  p_d \quad \text{on } \Gm_{sd}, 
\label{normal.stress}\\
- (\bT(\bu_s,p_s) \bn_s) \cdot \btau_l \equiv
- 2 \nu_s (\bD(\bu_s)\bn_s)\cdot \btau_l
      & =  
\frac{\nu_s\alpha}{\sqrt{K_l}}
\bu_s\cdot\btau_l, \,\, 
1\leq l\leq d-1 \text{ on } \Gm_{sd}.
\label{BJS}
 \end{align}
Conditions \eqref{cont.flux} and \eqref{normal.stress} denote the
continuity of flux and normal stress through $\Gm_{sd}$,
respectively. Condition \eqref{BJS} is the well-known
Beaver-Joseph-Saffman slip with friction law \cite{BeaverJoseph, Saffman}, where
$\{\btau_l\}_{l=1}^{d-1}$ is an orthogonal system of unit tangent
vectors on $\Gm_{sd}$ and $K_l = (K(\bx,\omega) \btau_l)\cdot \btau_l$. 
The constant $\alpha>0$ is a friction coefficient determined experimentally.

\subsection{Multi-region Karhunen-Lo\`eve (KL) expansion} \label{KLsec}
Let $Y(\bx,\omega)=\ln(K(\bx,\omega))$ be the log permeability and define 
$$
Y'(\bx,\omega):=Y(\bx,\omega)-E[Y](\bx).
$$
We assume that the Darcy region $D_d$ is a union of several 
non-overlapping KL regions $D_d =
\bigcup_{i=1}^{\nomega}D_{KL}^{(i)}$, where the stochastic structure
of every region is independent from the others. In other words, the
covariance between any pair of points from different KL regions is
zero. The permeability in each region is parameterized by a 
local KL expansion. The global stochastic space is decomposed 
correspondingly by
$$
\Om=\bigotimes_{i=1}^{\nomega}\Om^{(i)}.
$$
Therefore, for each event $\omega\in\Om$, we can write
$\omega=(\omega^{(1)},\ldots,\omega^{(\nomega)})$ and
$$
Y'(\bx,\omega)=\sum_{i=1}^{\nomega}Y^{(i)}(\bx,\omega^{(i)}).
$$
Each $Y^{(i)}$ has a physical support in $D_{KL}^{(i)}$ with a covariance function
\begin{equation*}
C_{Y^{(i)}}(\mathbf{x},\mathbf{\bar{x}})=E[Y^{(i)}(\mathbf{x},\omega^{(i)})Y^{(i)}(\mathbf{\bar{x}},\omega^{(i)})].
\end{equation*}
Since the covariance function is symmetric and positive definite, 
it can be decomposed into the series expansion
\begin{equation*}
C_{Y^{(i)}}(\mathbf{x},\mathbf{\bar{x}}) = \sum_{j=1}^{\infty}\lambda_j^{(i)} f_j^{(i)}(\mathbf{x})f_j^{(i)}(\mathbf{\bar{x}}),
\end{equation*}
where the eigenvalues and eigenfunctions $\lambda_j^{(i)}$, $f_j^{(i)}$ respectively, are computed by using $C_{Y^{(i)}}$ as the kernel of Fredholm integral equation
\begin{equation} \label{Fred}
\int_{D_{KL}^{(i)}} C_{Y^{(i)}}(\mathbf{x},\mathbf{\bar{x}}) 
f_j^{(i)}(\mathbf{x}) d\mathbf{x} = 
\lambda_j^{(i)} f_j^{(i)}(\mathbf{\bar{x}}).
\end{equation}
The eigenfunctions of $C_{Y^{(i)}}$ are mutually orthogonal and form a complete spanning
set. Therefore the Karhunen-Lo\`eve expansion for the log permeability can be expressed as
\begin{equation} \label{KL}
Y'(\mathbf{x},\omega) = \sum_{i=1}^{\nomega} ~ \sum_{j=1}^{\infty} ~ \xi_j^{(i)}(\omega^{(i)}) \sqrt{\lambda_j^{(i)}} f_j^{(i)}(\mathbf{x}),
\end{equation}
where the eigenfunctions $f_j^{(i)}(\mathbf{x})$ computed in
\eqref{Fred} have been extended by zero outside of $D_{KL}^{(i)}$ and
$\xi_j^{(i)}:\Omega_i \to \mathbb{R}$ are independent identically
distributed random variables \cite{ghanem2003sf}. We assume that
$Y^{(i)}$ is Gaussian, so each $\xi_j^{(i)}$ is a
random variable with zero mean and unit variance, with a
probability density function $\rho_j^{(i)}(y) = \frac{1}{\sqrt{2 \pi}}
\exp[\frac{-y^2}{2}]$.

Since the
eigenvalues $\lambda_j^{(i)}$ decay rapidly with $j$ 
\cite{zhang2004eho}, it is reasonable to truncate the local KL expansions.
If the expansion is truncated prematurely, the
permeability may appear too smooth in a particular KL region. In our
case, for any KL region $i$, we truncate the expansion after its first
$\nterm(i)$ terms. Increasing $\nterm(i)$ introduces more
heterogeneity into the permeability realizations for a chosen
region. It is beyond the scope of this paper to address the modeling
error associated with truncating the KL expansion. Some work has been
done to quantify the modeling error \cite{Nob-Temp-Web} and it can be
reduced a posteriori \cite{chang2009comparative}. The truncation 
allows us to approximate \eqref{KL} by
\begin{equation} \label{KLapprox}
Y'(\mathbf{x},\omega) \approx \sum_{i=1}^{\nomega} \sum_{j=1}^{\nterm(i)}\xi_j^{(i)}(\omega^{(i)}) \sqrt{\lambda_j^{(i)}} f_j^{(i)}(\mathbf{x}).
\end{equation}
The above also shows that globally we have 
$\nterm := \sum_{i=1}^{\nomega} \nterm(i)$ terms in $Y'$. 

Let $\S_j^{(i)} = \xi_j^{(i)}(\Omega^{(i)})$ be the images of the
random variables. They form finite dimensional spaces that are local
to each KL region, $\S^{(i)} = \prod_{j=1}^{\nterm(i)}\S_j^{(i)}
\subseteq \mathbb{R}^{\nterm(i)}$, as well as a global vector space 
$\S = \prod_{i=1}^{\nomega}~\S^{(i)} \subseteq \mathbb{R}^{\nterm}$.
Let $\kappa$ be a function that provides a natural ordering
for the global number of stochastic dimensions.  Then the $j$-th
stochastic parameter of the $i$-th KL region have a global index in
$\{1,\ldots,\nterm\}$ given by the function
\begin{equation*}
\kappa(i,j) = \left\{ 
\begin{tabular}{ll}
$j$, & if $i=1$, \\
$j + \displaystyle\sum_{k=1}^{i-1} \nterm(k)$, & if $i>1$.
\end{tabular}
\right.
\end{equation*}
The joint PDF for $\mathbf{\xi} = (\xi_j^{(i)})_{\kappa}$ is $\rho = \prod_i \prod_j
\rho_j^{(i)}$.  According to \eqref{KLapprox}, we approximate
$Y(\mathbf{x},\omega)$ by $Y(\mathbf{x},\mathbf{y})$, where
$\mathbf{y} = \left(\xi_j^{(i)}(\omega^{(i)})\right)_{\kappa}$. We
further note that, since in most cases closed-form eigenfunctions and
eigenvalues are not readily available, 
the integral equation \eqref{Fred} is solved numerically.

\section{Multiscale mortar mixed finite element method}\label{MMMFEM}
\subsection{Domain decomposition}

The coupled Stokes-Darcy flow problem is solved using domain decomposition
method following the approach described in \cite{VasWangYot}. The
Stokes and Darcy domains are partitioned into $N_s$ and $N_d$
non-overlapping subdomains, respectively. Let $N=N_s+N_d$, with $D_s =
\bigcup_{i=1}^{N_s} D_i$, $D_d = \bigcup_{i=N_s+1}^{N} D_i$, and $D_i
\cap D_j = \emptyset$ for $i \neq j$. Let the interface between
adjacent subdomains be $\Gamma_{i,j} = \partial D_i \cap \partial
D_j$. Depending on the models of adjacent domains, we group all
interfaces into three different types: Stokes-Stokes type, Darcy-Darcy
type and Stokes-Darcy type, denoted by $\Gamma_{ss}$, $\Gamma_{dd}$
and $\Gamma_{sd}$, respectively. The union of all interfaces is 
defined as $\Gamma = \Gamma_{ss}\bigcup \Gamma_{dd} \bigcup
\Gamma_{sd}$. Let $\Gamma_i = \d D_i \cap \Gamma$.
In addition, it is assumed that in $D_d$ each KL region
is an exact union of subdomains. Let
$(\bu_i,p_i)=(\bu_s|_{D_i},p_s|_{D_i})$ if $D_i$ is a Stokes subdomain
and $(\bu_i,p_i)=(\bu_d|_{D_i},p_d|_{D_i})$ if it is a Darcy
subdomain.

In addition to the interface conditions \eqref{cont.flux}--\eqref{BJS} on $\Gamma_{sd}$,
we impose continuity of velocity and normal stress on $\Gamma_{ss}$, and continuity of normal
velocity and pressure on $\Gamma_{dd}$:
\begin{align}
& [\bu] = 0, \quad  [\mathbf{T}(\bu,p)\bn] = 0 \quad \textrm{on}~\Gamma_{ss}, \label{intcond1} \\
& [\bu\cdot\bn]  = 0,\quad [p]  = 0 \quad  \textrm{on}~\Gamma_{dd}, \label{intcond2}
\end{align} 
where $[ \ \cdot\ ]$ represents the jump on the interface. For
example, on $\Gm_{ij}$, $[\bu]=(\bu_i-\bu_j)|_{\Gm_{ij}}$,
$[\bu\cdot\bn]=\bu_i\cdot\bn_i+\bu_j\cdot\bn_j$, with $\bn_i$ the outer unit normal vector 
to $\partial D_i$.



\subsection{Variational formulation}

We make use of the usual notation for Hilbert
spaces $H^k(G)$ for a set $G \subset \mathbb{R}^d$. The $L^2(G)$
inner product is denoted by $(\cdot,\cdot)_{G}$ for scalar, vector and
tensor valued functions. For a section of a subdomain boundary $S$ we
write $\langle \cdot, \cdot \rangle_{S}$ for the $L^2(S)$ inner
product (or duality pairing). 

In the deterministic setting, the velocity and pressure spaces are
defined as
\begin{align*}
\widetilde{V}_s &= \{\mathbf{v} \in L^2(D_s)^d : \mathbf{v}|_{D_{i}} \in H^1(D_{i})^d, 
1 \leq i \leq N_s,\  \mathbf{v} =  \mathbf{0}  \mbox{ on } \Gamma_s\}, \ \ \widetilde{W}_s = L^2(D_s),
\\
\widetilde{V}_d & = \{\mathbf{v} \in L^2(D_d)^d : \mathbf{v}|_{D_{i}} \in H(\mbox{div}; D_{i}), \,
N_s+1 \leq i \leq N, \, \mathbf{v}\cdot \mathbf{n}_d = 0  \mbox{ on } \Gamma_d \}, 
\ \ \widetilde{W}_d = L^2(D_d), \\
\widetilde{V} & = \widetilde{V}_s \oplus \widetilde{V}_d, \ \ 
\widetilde{W} = (\widetilde{W}_s \oplus \widetilde{W}_d) \cap L_0^2(D),
\end{align*}
where
$$
H({\rm div};D_{i}) = \{\bv_{d,i} \in L^2(D_{i})^d : \div \bv_{i} \in L^2(D_{i}) \}, 
$$
equipped with the norm
$$\|\bv\|_{H({\rm div};D_{i})}=(\|\bv\|_{L^2(D_{i})}^2+\|\div \bv\|_{L^2(D_{i})}^2)^{1/2},$$
and $L_0^2(D)$ is the space of functions with mean value zero on $D$. We also introduce 
a space for the Lagrange multiplier to impose continuity on the interfaces, 
$$
\widetilde{\Lambda} = \widetilde{\Lambda}_{ss}\oplus\widetilde{\Lambda}_{sd}
\oplus\widetilde{\Lambda}_{dd},
$$
$$
\widetilde{\Lambda}_{ss} = \left(\widetilde{V}_S|_{\Gm_{ss}}\right)', \quad
\widetilde{\Lambda}_{sd} = \left(\widetilde{V}_D\cdot\n|_{\Gm_{sd}}\right)', \quad
\widetilde{\Lambda}_{dd} = \left(\widetilde{V}_D\cdot\n|_{\Gm_{dd}}\right)', 
$$
where $(\cdot)'$ denotes the dual space. On $\Gamma_{ss}$ the Lagrange
multiplier has the physical meaning of normal stress and on
$\Gamma_{dd}\bigcup \Gamma_{sd}$ it has the meaning of pressure.

In the space of stochastic dimensions $\mathbb{S}$ we define
$$
L^2_{\rho}(\S) = \left\{ \bv:\S \to \R^d~\big| \left(\int_{\S}
\|\bv(\by)\|^2 \rhoydy\right)^{1/2} < \infty \right\},
$$
and take its tensor product with the deterministic spaces above:
$$
V = \widetilde{V}\otimes L^2_{\rho}(\S), 
\quad W = \widetilde{W} \otimes L^2_{\rho}(\S), \quad 
\Lambda = \widetilde{\Lambda} \otimes L^2_{\rho}(\S).
$$

The weak formulation for the coupled Stokes-Darcy problem
\eqref{stokes.1}--\eqref{BJS} and \eqref{intcond1}--\eqref{intcond2}
is given by: find $(\bu,p,\lambda)\in V\times W\times \Lambda$, such
that
\begin{align}\label{varform1}
a(\mathbf{u},\mathbf{v}) + b(\mathbf{v},p) + b_{\Lambda}(\mathbf{v},\lambda) 
&= \int_{\mathbb{S}} (\mathbf{f},\mathbf{v})_D \, \rho(\by)d\by, \quad \forall \,\mathbf{v} \in V,\\
b(\mathbf{u},w) &= -\int_{\mathbb{S}}(q_d,w)_{D_d} \, \rho(\by)d\by, \quad \forall \,w \in W,\\
b_{\Lambda}(\mathbf{u},\mu) &= 0, \quad \forall \,\mu \in \Lambda,
\label{varform2}
\end{align} 
where
\begin{align*}
&\tilde{a}_{i}(\mathbf{u}_{i},\mathbf{v}_{i}) 
 = 2\,\nu_s(\mathbf{D}(\mathbf{u}_{i}),\mathbf{D}(\mathbf{v}_{i}))_{D_i} 
 + \sum_{l=1}^{d-1} \left\<\frac{\nu_s\alpha}{\sqrt{K_l}}\mathbf{u}_{i}\cdot
 \btau_l,
 \mathbf{v}_i \cdot \btau_l\right\>_{\Gamma_i\cap\Gamma_{sd}},
 \, 1 \leq i \leq N_s, \\ 
&\tilde{a}_{i}(\mathbf{u}_{i},\mathbf{v}_{i})
 = \nu_d(K^{-1}\mathbf{u}_{i},\mathbf{v}_{i})_{D_i}, \ N_s+1 \leq i \leq N, \\
&  \tilde{b}_i(\mathbf{v}_i,w_i) 
 = - (\div \mathbf{v}_i,w_i)_{D_i}, \ 1 \leq i \leq N, \\
&\tilde{a}(\mathbf{u},\mathbf{v})  = \sum_{i=1}^{N}\tilde{a}_{i}(\mathbf{u}_i,\mathbf{v}_i), 
\quad 
\tilde{b}(\mathbf{v},w) = \sum_{i=1}^{N}\tilde{b}_i(\mathbf{v}_i,w_i),\\
&\tilde{b}_{\Lambda}(\mathbf{v},\mu)
 = \<[\mathbf{v}],\mu\>_{\Gamma_{ss}} + \<[\mathbf{v}\cdot\mathbf{n}],\mu\>_{\Gamma_{dd}}  + 
\<[\mathbf{v}\cdot\mathbf{n}],\mu\>_ {\Gamma_{sd}}, \\
& a(\mathbf{u},\mathbf{v}) 
 = \int_{\S}\tilde{a}(\mathbf{u},\mathbf{v})\rho(\by)d\by, 
\quad b(\bu,\bv)=\int_{\S}\tilde{b}(\mathbf{u},\mathbf{v})\rho(\by)d\by, \quad 
b_{\Lambda}(\mathbf{v},\mu) = \int_\mathbb{S}\tilde{b}_{\Lambda}(\bv,\mu)\rho(\by)d\by.
\end{align*}  
We note that in \eqref{varform2} the continuity of velocity or normal flux is imposed weakly
on $\Gamma\times\S$ via the Lagrange multiplier.

\subsection{Finite element discretization}

For the discretization of the stochastic variational formulation
\eqref{varform1}--\eqref{varform2}, we employ the multiscale mortar
mixed finite element method (MMMFEM) in the physical dimensions $D$,
coupled with the stochastic collocation method, such as tensor product
or sparse grid collocation using Gauss-Hermite quadrature rule for the
stochastic dimensions $\S$. Therefore, approximating
\eqref{varform1}--\eqref{varform2} involves solving a sequence of
independent deterministic problems using a domain decomposition
algorithm for the Stokes-Darcy coupled problem.  The resulting
solutions are realizations in stochastic space and weights in
quadrature rules for computing statistical moments.

The MMMFEM allows for non-conforming grids along the subdomain
interfaces. Each subdomain $D_i$ is covered by a shape regular finite
element partition $\mathcal{T}_{h_i}$ with $h_i$ being the maximal
element diameter. Let $\mathcal{T}_{h} = \bigcup_i \mathcal{T}_{h_i}$
be the global fine mesh and let $h=\max_{i=1}^N h_i$. Let
$\widetilde{V}_{h,i} \times \widetilde{W}_{h,i} \subset
\widetilde{V}_i \times \widetilde{W}_i$ be any stable
pair of finite element spaces on $D_i$ for Stokes, $1\leq i\leq N_s$, or
Darcy, $N_s+1\leq i\leq N$, satisfying the inf-sup condition
 \begin{equation}\label{inf-sup.S}
 \inf_{0 \neq w_{h,i} \in \widetilde{W}_{h,i}} \ \ \sup_{0 \neq \bv_{h,i} \in \widetilde{V}_{h,i}}
 \frac{( w_{h,i} , \div \bv_{h,i} )_{D_i}}{\|\bv_{h,i}\|_{H^1(D_i)} \ \ 
\|w_{h,i}\|_{L^2(D_i)}} \geq \beta_s > 0.
 \end{equation}
 where $\widetilde{V}_i = \widetilde{V}|_{D_i}$ and
 $\widetilde{W}_i = \widetilde{W}|_{D_i}$. Examples include the 
Taylor-Hood elements \cite{TH}, the MINI elements
 \cite{ABF} and the conforming Crouzeix-Raviart elements
 \cite{CR} for Stokes, as well as 
the Raviart-Thomas (RT) spaces \cite{RT}, and the Brezzi-Douglas-Marini (BDM) spaces
 \cite{BDM} for Darcy. The global discrete velocity and pressure spaces are
 given as $\widetilde{V}_h = \bigoplus_{i=1}^N \widetilde{V}_{h,i}$
 and $\widetilde{W}_h = \bigoplus_{i=1}^N \widetilde{W}_{h,i}$,
 respectively.
 
Next, each subdomain interface $\Gamma_{i,j}$ is partitioned into a
coarse $(d-1)$-dimensional quasi-uniform affine mesh
$\mathcal{T}_{H_{i,j}}$ with maximal element diameter $H_{i,j}$.  A
mortar space $\widetilde{\Lambda}_{H_{i,j}}(\Gm_{i,j})\subset
L^2(\Gm_{i,j})$ is defined to impose weakly the continuity of the discrete normal
flux or velocity across the non-matching grids. It may
contain continuous or discontinuous piecewise polynomials. Let $H=
\max_{i,j}H_{i,j}$ and $\widetilde{\Lambda}_H = \bigoplus_{i,j}
\widetilde{\Lambda}_{H_{i,j}}$. The semidiscrete approximation of
\eqref{varform1}--\eqref{varform2} is: find
$$
\bu_h: \S \to \widetilde{V}_h , \quad
p_h:\S \to \widetilde{W}_h, \quad
\lambda_H: \S \to \widetilde{\Lambda}_H 
$$
such that for $\rho$-almost every $\by\in\S$, the following deterministic problem holds:
\begin{align}\label{semidisc1}
\tilde{a}(\bu_h,\bv_h) + \tilde{b}(\bv_h,p_h) + \tilde{b}_{\Lambda}(\bv_h,\lambda_H) 
&=(\mathbf{f},\bv_h)_D, \quad \forall \,\bv_h \in \widetilde{V}_h,\\
\tilde{b}(\bu_h,w_h) &= -(q_d,w_h)_{D_d}, \quad \forall \,w_h \in \widetilde{W}_h,\\
\tilde{b}_{\Lambda}(\bu_h,\mu) &= 0, \quad \forall \,\mu \in \widetilde{\Lambda}_H.
\label{semidisc2}
\end{align} 
For the well-posedness of \eqref{semidisc1}--\eqref{semidisc2}, we
refer the reader to \cite{LaytonSchieweckYotov,GVY}. The mortar
variable $\lambda_H$ represents the stress on $\Gm_{ss}$ and pressure
on $\Gm_{dd}$. Therefore the continuity of normal stress and pressure
in \eqref{intcond1}--\eqref{intcond2} is satisfied on the mortar
mesh. On the other hand, the continuity of velocity on $\Gm_{ss}$ and
normal flux on $\Gm_{dd}$ in \eqref{intcond1}--\eqref{intcond2} is
satisfied weakly on the coarse scale $H$ by \eqref{semidisc2}. The use
of a coarse Lagrange multiplier results in a less expensive interface
problem and gain in computational efficiency.

\section{Stochastic collocation}\label{collocation}

The stochastic collocation method is employed to approximate the semidiscrete
solution $(\bu_h,p_h,\lambda_H)$ by an interpolant $\mathcal{I}_m$,
where $m$ is a multi-index indicating the desired
polynomial degree of accuracy in stochastic dimensions. It is uniquely
formed on a set of $\nreal$ stochastic points ${\by_k}$, where
$\nreal$ is a function of $m$ and the fully discrete solution is given
by
$$
\bu_{h,m}(\mathbf{x},\mathbf{y}) = \mathcal{I}_m \bu_h(\mathbf{x},\mathbf{y}), 
\quad
p_{h,m}(\mathbf{x},\mathbf{y}) = \mathcal{I}_m p_h(\mathbf{x},\mathbf{y}), 
\quad 
\lambda_{H,m}(\mathbf{x},\mathbf{y}) = \mathcal{I}_m \lambda_H(\mathbf{x},\mathbf{y}).
$$
Consider the Lagrange basis $\{L_m^{\{k\}}(\mathbf{y})\}$, which
satisfies $\{L_m^{\{k\}}(\mathbf{y}_j)\} = \delta_{kj}$. The Lagrange
representation for the fully discrete solution is
\begin{equation}\label{Lagr}
(\bu_{h,m}, p_{h,m}, \lambda_{H,m})(\mathbf{x},\mathbf{y}) = 
\sum_{k=1}^{\nreal} (\bu_h^{\{k\}}, p_h^{\{k\}}, \lambda_H^{\{k\}})
(\bx) L_m^{\{k\}}(\mathbf{y}),
\end{equation}
where $(\bu_h^{\{k\}},p_h^{\{k\}},\lambda_H^{\{k\}})$ is the evaluation
of semidiscrete solution $(\bu_h,p_h,\lambda_H)$ at the point in
stochastic space $\mathbf{y}_k$. Hence, computing \eqref{Lagr} requires
solving a deterministic problem for each permeability 
realization $K^{\{k\}}(\bx) = K(\bx,\by_k)$, $k=1,\ldots,\nreal$:
find 
$\bu_h^{\{k\}} \in \widetilde{V}_h$, 
$p_h^{\{k\}} \in \widetilde{W}_h$, and
$\lambda_H^{\{k\}}\in \widetilde{\Lambda}_H $ such that
\begin{align}\label{fullydisc1}
\tilde{a}(\bu_h^{\{k\}},\bv_h) + \tilde{b}(\bv_h,p_h^{\{k\}}) 
+ \tilde{b}_{\Lambda}(\bv_h,\lambda_H^{\{k\}}) 
&=(\mathbf{f},\bv_h)_D, \quad \forall \,\bv_h \in \widetilde{V}_h,\\
\tilde{b}(\bu_h^{\{k\}},w_h) & = -(q_d,w_h)_{D_d}, \quad \forall \,w_h \in \widetilde{W}_h,\\
\tilde{b}_{\Lambda}(\bu_h^{\{k\}},\mu) &= 0, \quad \forall \,\mu \in \widetilde{\Lambda}_H.
\label{fullydisc2}
\end{align} 
We note that substituting the Lagrange representation \eqref{Lagr} into the 
expectation integral \eqref{exp} results in a quadrature rule for computing
the stochastic mean of the discrete solution with weights 
$w_m^{\{k\}} = \int_{\S} L_m^{\{k\}}(\by)\rhoydy$.
Different collocation methods could be obtained by
choosing different collocation points $\by_k$. We
consider tensor product and sparse grids, both of which are
constructed using one-dimensional rules, where the points in
dimension $\S^{(i)}_j$ are the zeros of the orthogonal polynomials with
respect to the $L^2_\rho(\S^{(i)}_j)$ inner product. In our case
the random variables are Gaussian and we 
choose the zeros of the $N(0,1)$ Hermite
polynomials
$$H_{m}(y) = m! \sum_{k=0}^{[m/2]}(-1)^k\frac{(2y)^{m-2k}}{k! (m-2k)!}.$$
The one dimensional Gauss-Hermite
quadrature rule is accurate to degree $2m-1$.  
Let the weights and abscissae for the quadrature rule be
$$
{\cal W}(m) = \{w^1_m,\ldots,w^m_m\} \quad\mbox{and}\quad
{\cal H}(m) = \{h^1_m,\ldots,h^m_m\}.
$$
These can be computed with a symbolic manipulation software package.

\subsection{Collocation on tensor product grids}
In tensor product collocation \cite{babuska2008stochastic},
the polynomial accuracy is prescribed
in terms of \emph{component} degree.  This allows for anisotropic
rules to be easily constructed, but the number of collocation points
grows exponentially with the number of dimensions or polynomial
accuracy.  Thus, tensor product collocation is usually used in
problems with relatively low number of stochastic dimensions.

Recall that $\nterm(i)$ is the stochastic dimension in KL region $i$ and
$\nterm$ is the total stochastic dimension. We choose $\ncoll(i,j)$
collocation points in stochastic dimension $j$ of KL region $i$,
and define $\bm=(\ncoll(i,j))_\kappa$, which is an
$\nterm$-dimensional multi-index, as the required component degree of
the interpolant in the stochastic space $\S$.  The corresponding
anisotropic tensor product Gauss-Hermite interpolant in
$\nterm$-dimensions is given by
\begin{eqnarray*}
\nonumber \mathcal{I}^{\textrm{TG}}_{\bm}f(\by) \!\! &=& \!\! (\mathcal{I}_{\bm(1)} \otimes \cdots \otimes \mathcal{I}_{\bm(\nterm)})f(\by) \\
&=& \!\! \sum_{k_1=1}^{\bm(1)}\cdots \!\!\! \sum_{k_{\nterm}=1}^{\bm(\nterm)} \!\!\! f(h_{\bm(1)}^{k_1},\ldots,h_{\bm(\nterm)}^{k_{\nterm}}) L_{\bm(1)}^{k_1}(y_1) \cdots L_{\bm(\nterm)}^{k_{\nterm}}(y_{\nterm}),
\end{eqnarray*}
which interpolates the semi-discrete solution into the polynomial
space $\mathbb{P}_{\bm} = \prod_k \mathbb{P}_{\bm(k)}$ in $\S$.
The set of abscissae for this rule is
\begin{equation*} 
\mathcal{T}(\bm) = \bigotimes_{k=1}^{\nterm} \mathcal{H}(\mathbf{m}(k)) = \bigotimes_{i=1}^{\nomega} \left(\bigotimes_{j=1}^{\nterm(i)} \mathcal{H}(\ncoll(i,j))\right),
\end{equation*}
and the weight for the point $(h_{\bm(1)}^{k_1},\ldots,h_{\bm(\nterm)}^{k_{\nterm}})$
is 
$
w(\mathbf{k})=\prod_{i=1}^{\nterm}w_{\bm(i)}^{k_i}
$.

\begin{figure}
\centerline{
\includegraphics[width=3in]{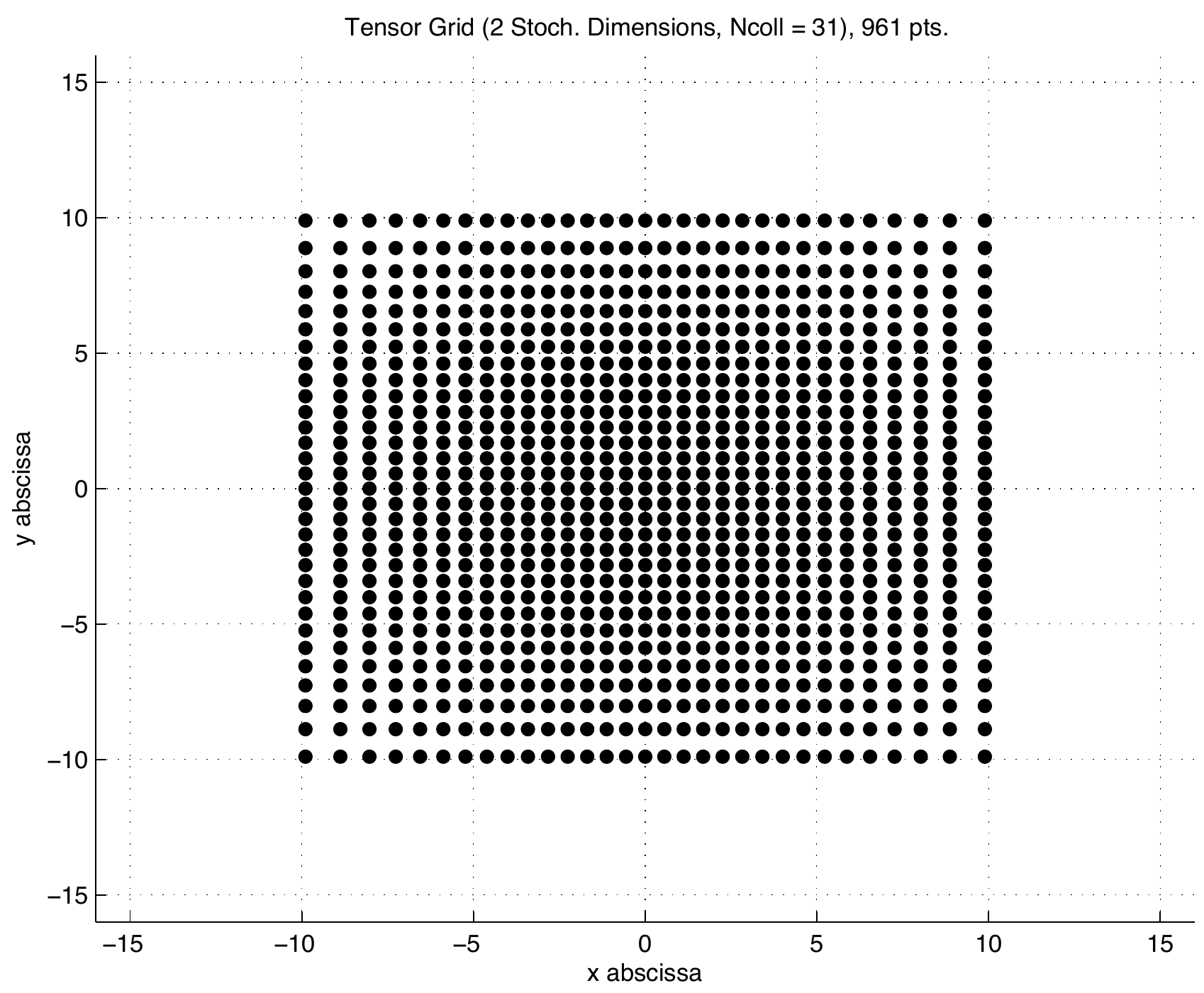}
\includegraphics[width=3in]{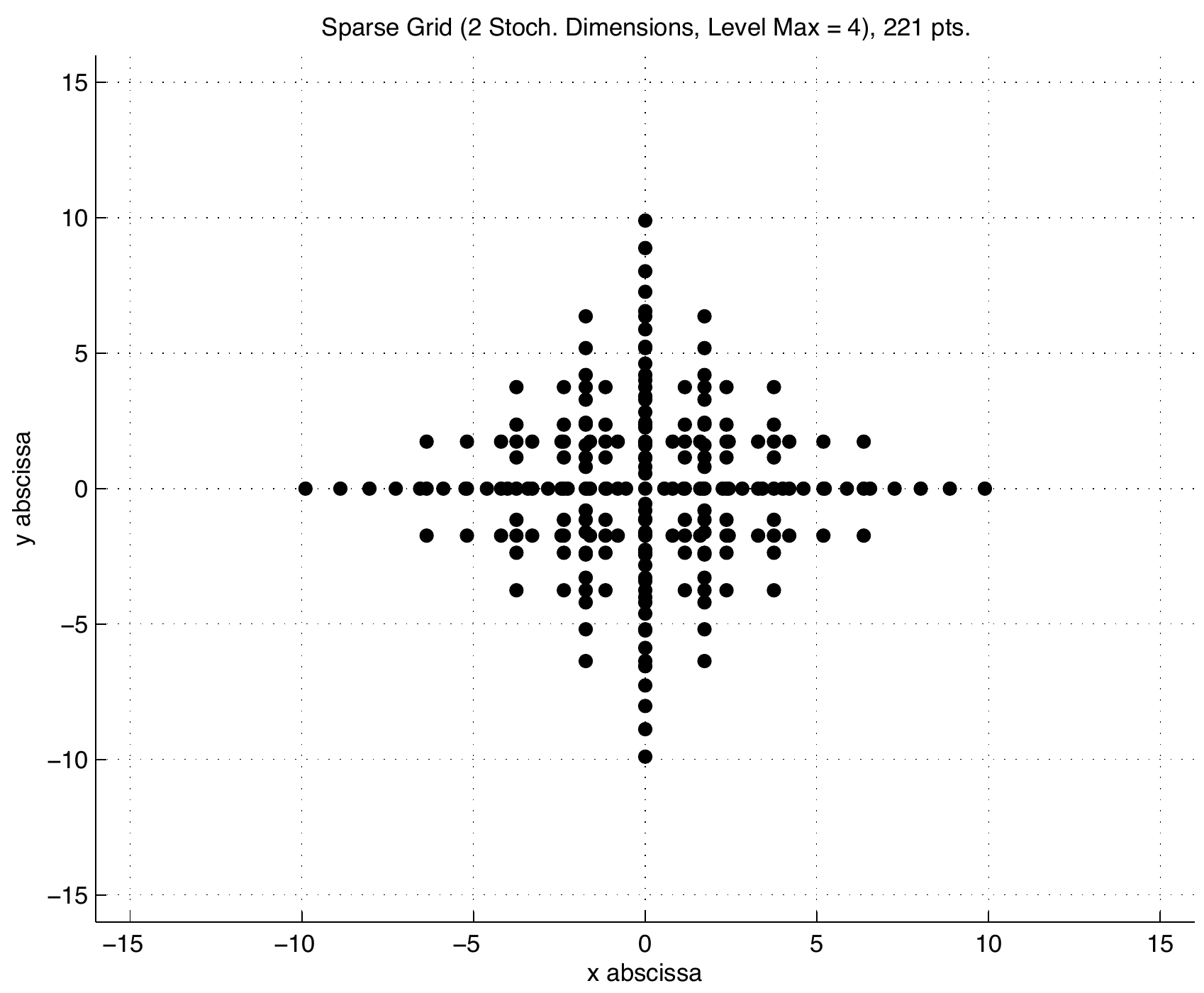}
}
\caption{A Gauss-Hermite tensor-product grid (left) versus a Gauss-Hermite sparse grid (right) with a comparable number of points on each axis.}
\label{grids}
\end{figure}

\subsection{Collocation on sparse grids}

In sparse grid collocation \cite{Xiu-Hest,Nob-Temp-Web}, the
polynomial accuracy is prescribed in terms of \emph{total} degree.
Sparse grids rules require much fewer points than tensor product rules
as the dimension increases, but have the same asymptotic
accuracy. This makes them applicable for problems with high number of
stochastic dimensions.  A tensor grid and a sparse grid rules that are
comparable in terms of accuracy are shown in Figure \ref{grids}.

Sparse grid rules are constructed hierarchically as linear combination
of tensor products of nested one-dimensional rules, so that the total
polynomial degree is a constant independent of dimension.  The
$\nterm$-dimensional sparse grid quadrature rule of level $\lmax$ is
accurate to degree $(2 \cdot \lmax+1)$.

Each level between $\lmax$ and $\lmin = \max \{0,\lmax-\nterm+1\}$ is
associated with a multi-index $\bp=(p_1,\ldots,p_{\nterm})$, where
\mbox{$|\bp| = \sum p_i$} denotes the levels of one dimensional rules
used in each stochastic dimension.  We consider the Gauss-Hermite
points $\mathcal{H}(2^{p_i+1}-1)$ as the one dimensional abscissae of
level $p_i$.  Level $0$ consists of a single point, and for every
subsequent level the number of points doubles plus one.
For each partition $\bp$, let the multi-index 
$\bm = 2^{\bp+\mathbf{1}}-\mathbf{1}$ denote its degree.
The $\nterm$-dimensional sparse grid Gauss-Hermite interpolant is given by
\begin{equation*}
\mathcal{I}^{\textrm{SG}}_{\lmax}f(\by) = \sum_{\lmin \leq |\bp| \leq \lmax} (-1)^{\lmax-|\bp|} \cdot \binom{\nterm-1}{\lmax-|\bp|} \cdot \mathcal{I}^{\textrm{TG}}_{\bm}f(\by),
\end{equation*}
and the set of abscissae is
\begin{equation*} 
\mathcal{S}(\lmin,\lmax,\nterm) = \bigcup_{\lmin \leq |\bp| \leq \lmax} \bigotimes_{i=1}^{\nterm} \mathcal{H}(2^{p_i+1}-1).
\end{equation*}

\section{Collocation-MMMFEM algorithms for Stokes-Darcy}\label{algorithms}

As noted above, computing the fully discrete approximation to
\eqref{varform1}--\eqref{varform2} in the coupled physical-stochastic
space requires solving a sequence of deterministic Stokes-Darcy
problems \eqref{fullydisc1}--\eqref{fullydisc2}. These are solved in
parallel by using a non-overlapping domain decomposition algorithm
developed in \cite{VasWangYot}. It reduces to global fine scale
problem to a coarse scale symmetric and positive definite interface
problem, which is solved by a Krylov iterative solver such as the
Conjugate Gradient (CG) method. We develop three algorithms that
combine different implementations of the domain decomposition
algorithm with stochastic collocation. The first method is based on
the traditional domain decomposition algorithm, which requires solving
subdomain problems at every CG iteration. The number of subdomain
solves, which is the dominant cost in the algorithm, is proportional
to the product of the total number of stochastic collocation points
and the number of interface iterations.  The other two methods exploit
the efficiency of precomputing a multiscale flux basis, which is used
to replace the subdomain solves by linear combination of the
basis functions. One of these methods, referred to as deterministic multiscale
flux basis, computes a new multiscale basis at each stochastic
collocation. The number of subdomain solves for this method is
proportional to the product of the total number of stochastic
collocation points and the number of mortar degrees of freedom per
subdomain.  For the last method we follow the approach in
\cite{GanisYotovZhong} to gain additional efficiency by exploding the
locality of the KL expansions. In particular, we precompute a full
multiscale flux basis on each subdomain for all local stochastic
realizations before the global collocation loop begins. This
stochastic basis is reused multiple times during the global
collocation loop, since the same local subdomain stochastic
realization occurs over multiple global realizations. In this method,
referred to as stochastic multiscale flux basis, the number of
subdomain solves is proportional to the product of the number of
mortar degrees of freedom per subdomain and the number of local
stochastic collocation points, resulting in additional computational
savings.

\subsection{Reduction to an interface problem} \label{reduction}
Before presenting the three algorithms, we describe the
non-overlapping domain decomposition algorithm for solving the
deterministic Stokes-Darcy problem
\eqref{fullydisc1}--\eqref{fullydisc2} that reduces the global problem
to a coarse scale interface problem. We
split \eqref{fullydisc1}--\eqref{fullydisc2} into two sets of
complementary subproblems, one with specified interface data and all 
other data set to zero, and the other with zero interface data and the 
given source term and outside boundary conditions.
For Stokes subdomains $D_{i}$, $1\leq i\leq N_s$, one of
the subproblems is to find
$\bu_{h,i}^{*,\k}(\bl) \in \widetilde{V}_{h,i}\slash \ker \,\tilde{a}_i$ and 
$p_{h,i}^{*,\k}(\bl) \in \widetilde{W}_{h,i}$
with specified $\bl=(\lambda_n,\bl_\tau)$, where $\lambda_n$ and
$\bl_\tau=(\lambda^1_{\tau},\ldots,\lambda^{d-1}_{\tau})$ represent
the normal stress and tangential stress on $\Gm_{ss}$, respectively,
such that 
\begin{align}
  & \tilde{a}_{i}(\bu_{h,i}^{*,\k}(\bl),\bv_i)
  + \tilde{b}_i(\bv_i,p_{h,i}^{*,\k}(\bl))
\nonumber \\
& \qquad\qquad = -\<\lambda_n,\bv_i\cdot\bn_i\>_{\Gamma_i}
- \sum_{l=1}^{d-1}\< \lambda_{\tau}^l , \bv_i\cdot\btau_{i}^l \>_{\Gamma_i\cap\Gm_{ss}},
\quad \forall \, \bv_i \in \widetilde{V}_{h,i}\slash \ker \,\tilde{a}_i,
\label{starSproblem1} \\
& \tilde{b}_i(\bu_{h,i}^{*,\k}(\bl),w_i) = 0,
\quad \forall \,w_i \in \widetilde{W}_{h,i},\label{starSproblem2}
\end{align}
where $\{\btau_{i}^l\}_{l=1}^{d-1}$ is an orthogonal set of unit
vectors tangential to $\d D_{i}$ and the kernel space $\ker\,
\tilde{a}_i := \{\bv \in \widetilde{V}_i: \tilde{a}_i(\bv,\bv) = 0\}$
consists of a subset of all rigid body motions depending on the
boundary types of $D_{i}$. For more detail we refer the reader to
\cite{VasWangYot}. The complementary subproblem is to find  
$\bar{\bu}^{\k}_{h,i} \in \widetilde{V}_{h,i}\slash\ker\,a_i$ and 
$\bar{p}^{\k}_{h,i} \in \widetilde{W}_{h,i}$ such that
\begin{align}\label{barSproblem1}
  & \tilde{a}_{i}(\bar{\bu}^{\k}_{h,i},\bv_i)
  + \tilde{b}_i(\bv_i,\bar{p}^{\k}_{h,i}) = (\ff_{i},\bv_i)_{D_i} ,
  \quad \forall \, \bv_i \in \widetilde{V}_{h,i}\slash\ker\,a_i,\\
& \tilde{b}_i(\bar{\bu}^{\k}_{h,i},w_i) = 0, \quad \forall \, w_i \in \widetilde{W}_{h,i}.
\label{barSproblem2}
\end{align}
Note that the first problem \eqref{starSproblem1}--\eqref{starSproblem2}
has specified normal stress on the interface with zero boundary condition and
source term, while the second problem
\eqref{barSproblem1}--\eqref{barSproblem2} has specified boundary
condition and source term and zero normal stress on the interface. Similarly,
on the Darcy subdomains $D_{i}$, $N_s+1 \leq i \leq N$, the first
subproblem is to find
$\bu_{h,i}^{*,\k}(\lambda) \in \widetilde{V}_{h,i}$ and 
$p_{h,i}^{*,\k}(\lambda) \in \widetilde{W}_{h,i}$
with specified pressure $\lambda$ such that
\begin{align}\label{starDproblem1}
&\tilde{a}_{d,i}(\bu_{h,i}^{*,\k}(\lambda),\bv_i) + \tilde{b}_i(\bv_i,p_{h,i}^{*,\k}(\lambda))
  =-\<\lambda,\bv_i\cdot\bn_i\>_{\Gamma_i},
  \quad \forall \, \bv_i \in \widetilde{V}_{h,i},\\
  &\tilde{b}_i(\bu_{h,i}^{*,\k}(\lambda),w_i) = 0,
  \quad \forall \,w_i \in \widetilde{W}_{h,i}.\label{starDproblem2}
\end{align}
The complementary subproblem is to find 
$\bar{\bu}^{\k}_{h,i} \in \widetilde{V}_{h,i}$ and
$\bar{p}^{\k}_{h,i} \in \widetilde{W}_{h,i}$ such that
\begin{align}\label{barDproblem1}
  & \tilde{a}_{i}(\bar{\bu}^{\k}_{h,i},\bv_i)
  + \tilde{b}_i(\bv_i,\bar{p}^{\k}_{h,i}) = (\ff_{i},\bv_i)_{D_i} ,
  \quad \forall \, \bv_i \in \widetilde{V}_{h,i},\\
& \tilde{b}_i(\bar{\bu}^{\k}_{h,i},w_i) = -(q_d,w_i)_{D_i}, \quad \forall \, w_i \in \widetilde{W}_{h,i}.
\label{barDproblem2}
\end{align}
Using the interface condition \eqref{fullydisc2}, it is easy to see that
\eqref{fullydisc1}--\eqref{fullydisc2} is equivalent to solving
the interface problem: find $\lambda^{\k}_H \in \widetilde{\Lambda}_H$
such that 
\begin{equation}\label{interface}
s_H(\lambda^{\k}_H,\mu) := - \tilde{b}_\Lambda(\bu^{*,\k}_h(\lambda^{\k}_H),\mu) = 
\tilde{b}_\Lambda(\bar{\bu}^{\k}_h, \mu), \quad \forall \,
\mu \in \widetilde{\Lambda}_H,
\end{equation}
and recovering the global velocity and pressure by
$$
\bu^{\k}_h=\bu^{*,\k}_h(\lambda^{\k}_H)+\bar{\bu}^{\k}_h, \quad
p^{\k}_h=p_h(\lambda^{\k}_H)+\bar{p}^{\k}_h.
$$
It is shown in \cite{VasWangYot} that the bilinear form $s_H(\cdot,\cdot)$
in \eqref{interface} is symmetric and positive definite. We employ the CG
method for its solution.

The interface problem \eqref{interface} can be rewritten in an operator form
\begin{equation}\label{operatorform}
S_H^{\{k\}}\lambda_H^{\{k\}}=g_H^{\{k\}},
\end{equation}
where $S_H^{\{k\}}:\widetilde{\Lambda}_H \rightarrow \widetilde{\Lambda}_H$
is a Steklov--Poincar\'e type operator, defined as 
$$
\left(S_H^{\{k\}}\lambda,\mu\right)=s_H^{\{k\}}(\lambda,\mu),\quad
\forall \lambda,\mu\in \widetilde{\Lambda}_H,
$$
and $g_H^{\{k\}}$ is defined by
$\<g_H^{\{k\}},\mu\>_\Gamma = \tilde{b}_\Lambda(\bar{\bu}_h^{\{k\}}, \mu),
\ \ \forall\mu\in\widetilde{\Lambda}_H$.

It is important to note that the dominant cost in solving
\eqref{operatorform} is computing the action of operator $S_H^{\{k\}}$
in every CG iteration. This action involves solving the 
Neumann-to-Dirichlet problem
\eqref{starSproblem1}--\eqref{starSproblem2} in the Stokes subdomains,
and the Dirichlet-to-Neumann problem
\eqref{starDproblem1}--\eqref{starDproblem2} in the Darcy subdomains.

Let $\widetilde{\Lambda}_{H,i}$ be the restriction of
$\widetilde{\Lambda}_{H}$ to $\Gamma_i$ and let $S_{H,i}^{\{k\}}$
be the restriction of $S_{H}^{\{k\}}$ to $\widetilde{\Lambda}_{H,i}$. The
action of $S_{H}^{\{k\}}$ can be written as the sum of the actions of the
local operators:
$$
S_H^{\{k\}}\lambda_H = \sum_{i=1}^{N}S_{H,i}^{\{k\}}\lambda_{H,i}.
$$
We proceed with the description of the three methods for solving
\eqref{operatorform}.

\subsection{Method S1: Collocation without multiscale flux basis}

\begin{mdframed}
\textbf{Method S1} (without multiscale flux basis)

For $k = 1,\ldots, N_{real}$, do

\begin{tabular}{p{0.13\textwidth}p{0.7\textwidth}}
\quad Step 1: & Generate a permeability realization $K^{\{k\}}$ corresponding to the global collocation index $k$. \\
\quad Step 2: & Solve \eqref{operatorform} using traditional MMMFEM for $\lambda_H^{\{k\}}$.\\
\quad Step 3: & Add the solution to the statistical moments with the 
collocation weight applied.
\end{tabular}

End do 
\end{mdframed}

\vspace{2EX} In any subdomain, Step 2 costs one set of subdomain solves in
every CG iteration for applying $S_H^{\{k\}}$. Given a mortar function
$\lambda_H$, the action of $S_{H,i}^{\{k\}}\lambda_{H,i}$ in subdomain $D_i$
includes:
\begin{enumerate}
\item Project $\lambda_{H,i}$ onto subdomain boundaries:
  $\lambda_{h,i}=\mathcal{L}_{h,i}\lambda_{H,i}$,
  where $\mathcal{L}_{h,i}$ is the $L^2$-projection operator onto the
  velocity space (for Stokes) or the normal trace of the velocity space
  (for Darcy) on $\Gamma_i$. 
\item If $D_i$ is a Darcy subdomain, solve the subdomain problem
  \eqref{starDproblem1}--\eqref{starDproblem2} with Dirichlet data
  $\lambda_{h,i}$; If $D_i$ is a Stokes subdomain, solve the subdomain
  problem \eqref{starSproblem1}--\eqref{starSproblem2} with Neumann
  data $\lambda_{h,i}$.
\item Project the resulting flux in Darcy or velocity in Stokes back
  to the mortar space and compute the jump across the interfaces.
\end{enumerate}

Let $\cgiter(k)$ be the number of CG iterations for the $k$th
stochastic realization. Then in any subdomain $D_i$ the leading term
in the number of solves for method S1
is
\begin{equation}\label{S1-cost}
N_{S1}:=\sum_{k=1}^{\nreal} \cgiter(k),
\end{equation}
omitting the extra solves in the right hand side and solution
recovery.  The number of CG iterations depends on the condition number
of the interface operator $S_H^{\{k\}}$, which grows with refining the
grids or increasing the number of subdomains or permeability contrast
\cite{VasWangYot}, resulting in possibly large computational cost of
Method~S1.

\subsection{Method S2: Collocation with deterministic multiscale flux basis}
One way to reduce the cost of the MMMFEM for Stokes-Darcy flow is to
employ a multiscale flux basis following the method introduced in
\cite{GanVasWangYot}. When solving the deterministic problem
$\eqref{operatorform}$, before the CG iteration starts, we compute
on each subdomain the local solution for every mortar degree of freedom
associated with the subdomain and store the 
results as a basis of flux responses in Darcy or velocity responses
in Stokes. Then to evaluate the action of $S_H^{\{k\}}$ in every
interface iteration, we simply use the linear combination of the basis
instead of performing an extra subdomain solve. The algorithm for this
method is as follows:

\begin{mdframed}
\textbf{Method S2} (with deterministic multiscale flux basis)

For $k = 1,\ldots, N_{real}$, do

\begin{tabular}{p{0.13\textwidth}p{0.74\textwidth}}
\quad Step 1: & Generate a permeability realization $K^{\{k\}}$ corresponding to the global collocation index $k$. \\
\quad Step 2: & On each $D_i$, compute and save the multiscale flux basis
$\{\phi_{H,i}^{j,\{k\}}\}_{j=1}^{\ndofi}$\\ 
\quad Step 3: & Solve \eqref{operatorform} for  $\lambda_H^{\{k\}}$ using the
MMMFEM with the computed basis.\\
\quad Step 4: & Add the solution to the statistical moments with the
collocation weight applied.
\end{tabular}

End do 
\end{mdframed}

\vspace{2EX}

In Step 2, where $\ndofi$ is the number of mortar degrees of freedom on
$\Gamma_i$, we solve a subdomain problem for
each mortar basis function. The mortar basis
functions $\{\xi_{H,i}^{j}\}_{j=1}^{\ndofi}\subset
\widetilde{\Lambda}_{H,i}$ are defined such that any
$\lambda_{H,i}\in\widetilde{\Lambda}_{H,i}$ can be expressed uniquely
by their linear combination:
$$
\lambda_{H,i}=\sum_{j=1}^{\ndofi}c_i^j\xi_{H,i}^j.
$$   
For any $j=1\ldots\ndofi$, the multiscale flux basis is computed by
$$
\phi_{H,i}^{j,\{k\}}=S_{H,i}^{\{k\}}\xi_{H,i}^j.
$$
To evaluate $S_{H,i}^{\{k\}}\lambda_{H,i}$ in Step 3, we simply compute
$$
S_{H,i}^{\{k\}}\lambda_{H,i} =\sum_{j=1}^{\ndofi}c_i^j\phi_{H,i}^{j,\{k\}}.
$$   
Since there are no solves inside the CG loop, the number of
subdomain solves for
Method S2 does not depend on the number of CG iterations. In any
subdomain $D_i$, the number of subdomain solves for Method S2 has the leading
term
\begin{equation}\label{S2-cost}
  N_{S2}:=\ndofi*\nreal.
  \end{equation}
Unlike the cost of Method~S1 \eqref{S1-cost}, the dominant cost of
Method~S2 is independent of the condition number of the interface
operator, which makes it less sensitive to $h$, $H$, and jumps in $K$.
Comparing \eqref{S2-cost} to \eqref{S1-cost}, we observe that      
in every stochastic realization, Method S2 costs less than Method S1
if the number of CG iterations $\cgiter(k)$ is greater than the
maximum number of mortar degrees of freedom per subdomain.

\subsection{Method S3: Collocation with stochastic multiscale flux basis}

The third method, which utilizes the approach from
\cite{GanisYotovZhong}, is designed to achieve greater computational
savings than Method~S2 by exploiting the fact that the global KL
expansion is a sum of local KL expansions.  The main idea is to
compute and store a multiscale flux basis for all local KL
realizations before the global stochastic collocation loop.  Since a
local KL realization is repeated multiple times in the global KL
collocation loop, the precomputed stochastic multiscale flux basis can
be reused, resulting in further gain in computational efficiency.
Below we present the algorithm for method S3.

\begin{mdframed}
\textbf{Method S3} (with stochastic multiscale flux basis)

For any Darcy subdomain $D_i$ that belongs to KL region $j$,
set $\nreal(j)$ to be the number of
local KL realization; for any Stokes subdomain,
$\nreal(j)$ is set to 1.

\vspace{1.5EX}
\textbf{Pre-computation loop: stochastic multiscale basis}

For $k = 1,\ldots,\nreal(j)$, do

\begin{tabular}{p{0.13\textwidth}p{0.74\textwidth}}
\quad Step 1: & Generate a permeability realization $K^{\{k\}}$ corresponding to the local collocation index $k$. \\
\quad Step 2: & Compute and save the multiscale flux basis under the
current local index.\\ 
\end{tabular}

End do
\vspace{1.5EX}

\textbf{Main loop: global stochastic collocation}

For $l = 1,\ldots, N_{real}$, do

\begin{tabular}{p{0.13\textwidth}p{0.74\textwidth}}
\quad Step 3: & Generate a permeability realization $K^{\{l\}}$ corresponding to the global collocation index $l$. \\
\quad Step 4: & Convert the global index to the subdomain's local index $k$\\ 
\quad Step 5: & Solve \eqref{operatorform} for  $\lambda_H^{\{k\}}$ using
the MMMFEM with the computed basis for local index $k$ in step 2.\\
\quad Step 6: & Add the solution to the statistical moments with collocation
weight applied.
\end{tabular}

End do 
\end{mdframed}

\vspace{2EX}

The algorithms for index conversion in Step 4 for tensor product grid and
sparse grid are given in Algorithms 5 and 6 in
\cite{GanisYotovZhong}, respectively. Notice that all subdomain solves in
Method~S3,
except the extra two for RHS computation and solution recovery,
are done in the pre-computation loop only. Therefore in
subdomain $D_i$, the total number of solves has the leading term
\begin{equation}\label{S3-cost}
N_{S3}:=\ndofi*\nreal(j).
\end{equation}
Compared to the leading term in $N_{S2}$ for method S2
\eqref{S2-cost}, $N_{S3}$ is proportional to $\nreal(j)$ instead of
$\nreal$. When there exists more than one KL region in the problem
setting, in any KL region $j$ the number of local realizations
$\nreal(j)$ is less than the number of global realizations $\nreal$,
which makes Method~S3 cost less in terms of subdomain solves than
Method~S2.

\section{Numerical tests}\label{numerical}

In this section, three numerical examples are presented to test both
tensor product and sparse grid collocations with Methods S1, S2 and
S3. All three examples are motivated by coupled surface flow with
groundwater flow with uncertain non-stationary
heterogeneous permeability. The covariance
function in KL region $D_{KL}^{(i)}$ is
\begin{equation}\label{covariance}
C_{Y}^{(i)}(\bx,\overline{\bx})=\left(\sigma_Y^{(i)}\right)^2\exp\left[-\frac{|x_1-\overline{x}_1|}{\eta_1^{(i)}}-\frac{|x_2-\overline{x}_2|}{\eta_2^{(i)}}\right],
\end{equation}
where $\sigma_Y^{(i)}$ is the variance and $\eta_j^{(i)}$ is the
correlation length in the $j$th physical dimension. The three examples
are designed to test different levels of heterogeneity and uncertainty
in terms of number of KL regions, correlation length, and variance.
In addition, test case 2 illustrates the flexibility to use finer
physical discretization grids in KL regions with smaller correlation
lengths, and test case 3 presents an application to domains with
irregular geometry.

For all tests, we use Taylor-Hood \cite{TH} triangular finite element
in Stokes and the lowest order Raviart-Thomas \cite{RT} quadrilateral
finite elements in Darcy for space discretization. All interfaces are
discretized via discontinuous piecewise linear mortar finite
elements. The methods are implemented in parallel for distributed
memory computers using message passing interface (MPI). The numerical
tests are run on a parallel cluster of Intel
Xeon CPU E5-2650 v3 @ 2.30GHz processors with 192GB RAM, with each
subdomain assigned to one core.
 
\subsection{Test case 1}

\begin{figure}
\centering
\begin{minipage}[t]{.5\textwidth}
\centerline{\includegraphics[width=3.9in]{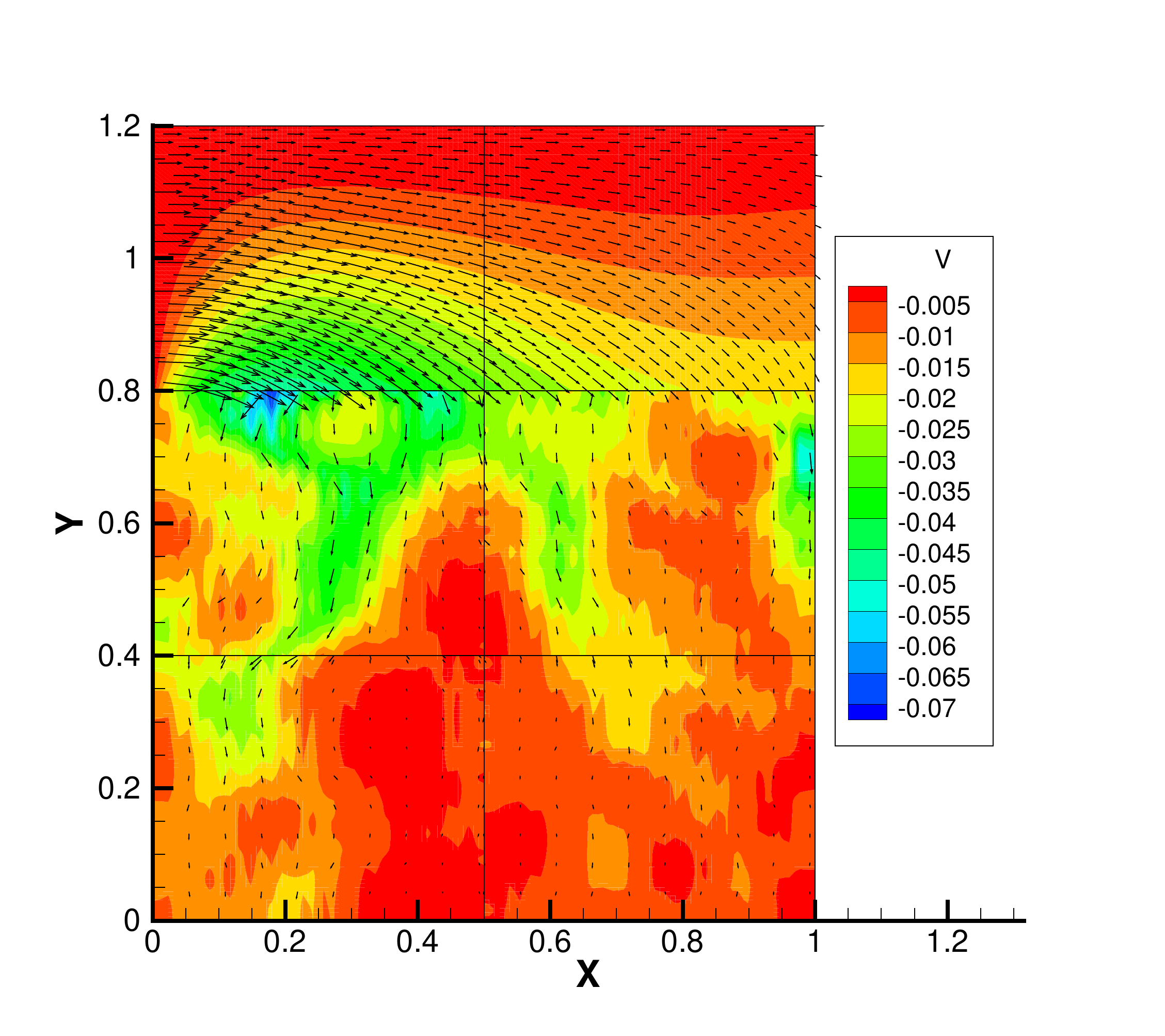}}
\end{minipage}%
\begin{minipage}[t]{.5\textwidth}
\centerline{\includegraphics[width=3.1in]{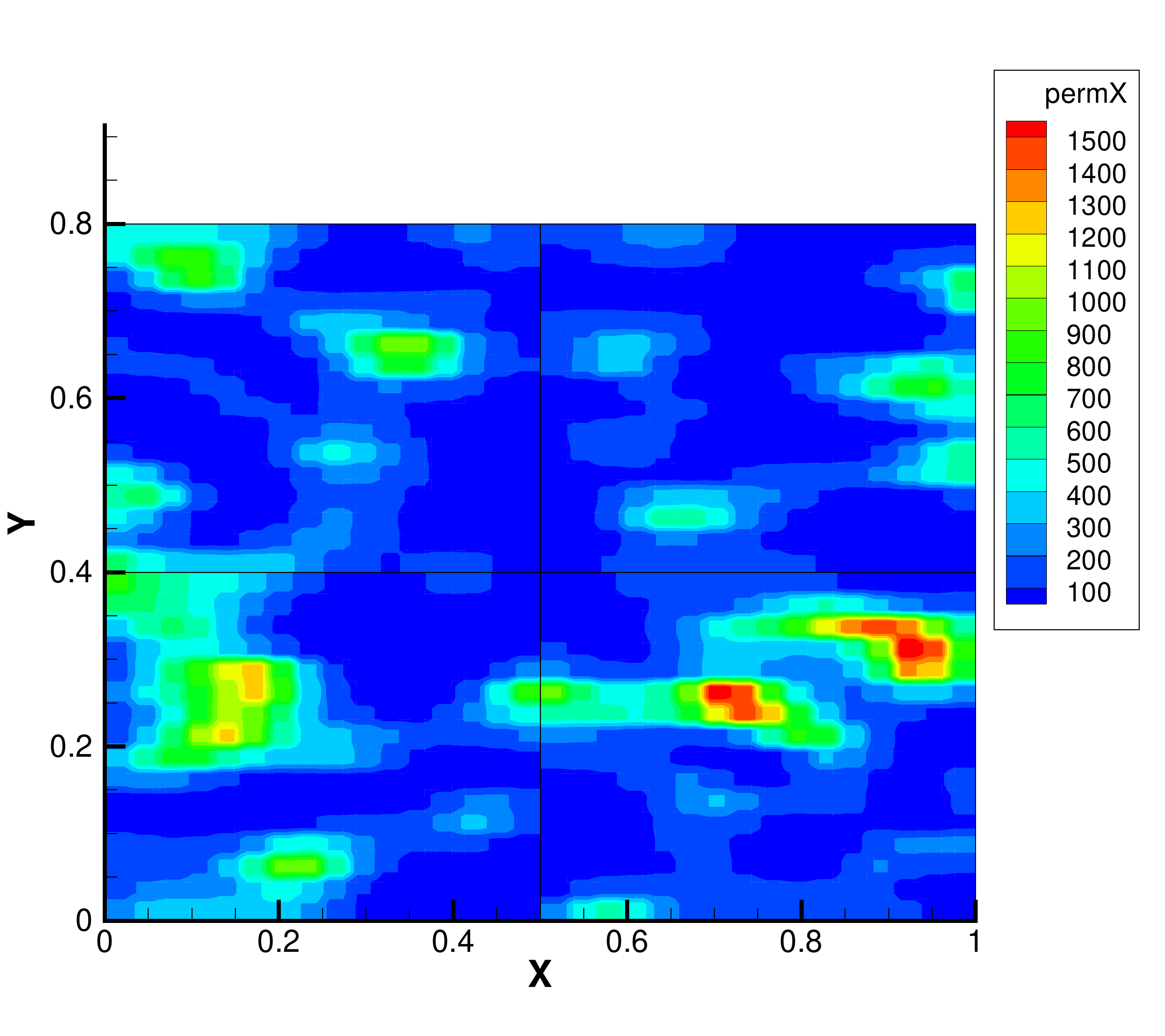}}
\end{minipage}
\caption{Case 1, solution realization - colors for vertical velocity and arrows for velocity vector (left) and permeability realization (right).}
\label{Case1perm}
\end{figure}

The first test case has a global domain $(0,1)\times(0,1.2)$, where
$D_d=(0,1)\times(0,0.8)$ and $D_s=(0,1)\times(0.8,1.2)$. The problem
is divided into 6 equal-size subdomains with two Stokes and four Darcy
subdomains, see Figure~\ref{Case1perm} (left).  The outside boundary
conditions are given as follows: in Darcy, zero pressure is specified
on the bottom edge, and no-flow condition on the left and right edges;
in Stokes, velocity is specified on the left and top edges, with a
parabolic inflow profile on the left and zero normal velocity and a
linearly decreasing tangential slip on the top, while zero normal and
tangential stresses are specified on the right edge.

The mean permeability is a given heterogeneous scalar field, which
varies two orders of magnitude, see Figure~\ref{Case1perm} (right).
There are two rectangular statistically independent KL regions in
Darcy, defined as $(0,1)\times(0,0.4)$ and $(0,1)\times(0.4,0.8)$. In
the tensor product grid collocation test we set $\nterm(1)=2\times 1 =
2$ and $\nterm(2)=3\times 3 = 9$. In the sparse grid test we set
$\nterm(1)=\nterm(2)=5\times 5=25$. In both KL regions we set
$\sigma_Y^{(i)} = 1.0$ and $\eta_j^{(i)} = 0.1$. In tensor product
collocation, the grid is isotropic with $\ncoll=2$ in each stochastic
dimension. In sparse grid collocation, $\lmax=1$, so on both cases we
have a third degree of accuracy in each stochastic dimension.
For the tensor product grid, the number of terms in the global KL expansion
is $\nterm = 2+9 = 11$, resulting in 
number of global realizations $\nreal= \ncoll^{\nterm} = 2^{11} = 2048$. For
the sparse grid, the number of terms in the global KL expansion
is $\nterm = 25+25 = 50$, resulting in number of global realizations
$\nreal= 2*\nterm + 1 = 101$.

For the physical space discretization, we have a local mesh of size $64
\times 64$ in every Darcy subdomain and size $32 \times 32$ in every Stokes
subdomain. The mortar mesh is $8 \times 1$ on every interface on
$\Gm_{dd}$ and $4\times 1$ on every interface on
$\Gm_{sd}\cup\Gm_{ss}$. One realization of the computed velocity field
is shown in Figure~\ref{Case1perm} left. We observe infiltration from
the Stokes region into the Darcy region with higher velocity in the
areas of higher permeability. The normal velocity across the
Stokes-Darcy interface is continuous and reasonable well resolved,
even with the fairly coarse mortar mesh. The tangential velocity along the
Stokes-Darcy interface is higher in the Stokes region.

\begin{table}
\begin{center}

{\em \textbf{ Tensor Product} Collocation, $\nterm=11$, $\ncoll=2$ (2048 realizations)
  
 }
\begin{tabular}{||l||c|c|c||}
\hline
& \textbf{Method S1} & \textbf{Method S2} & \textbf{Method S3} \\
\hline
Max. number of solves  & 1,499,257 & 167,936 & 45,056 \\ 
\hline
Runtime in seconds   & 33457.4   & 4464.9   & 3357.1 \\
\hline
\end{tabular}

\vspace{1.5ex}
 {\em \textbf{Sparse Grid} Collocation,  $\nterm=50$, $\lmax=1$ (101 realizations)
 
 }
\begin{tabular}{||l||c|c|c||}
\hline
& \textbf{Method S1} & \textbf{Method S2} & \textbf{Method S3}  \\
\hline
Max. number of solves  & 71,237 & 8,282 & 4,282  \\
\hline
Runtime in seconds   & 1566.1   & 220.5   & 156.8     \\ 
\hline
\end{tabular}
\end{center}
\caption{ Number of subdomain solves and runtime in seconds for Case 1. }
\label{Case1}
\end{table}

The number of solves and runtimes for all three methods on tensor
product and sparse grid collocations are shown in Table \ref{Case1}.
With tensor product collocation, Method S1 requires 1.5 million number
of solves. Recalling \eqref{S1-cost}, this very high number is due to
the high value of $\nreal = 2048$, as well as a large number of CG
iterations $\cgiter(k)$ in each realization.  In comparison, Methods
S2 requires only about 11\% of the number of solves and about 13\% of
the runtime of Method S1. We observe even larger savings in Method S3,
which requires only about 27\% of the number of solves of Method S2.
This is due the fact that the stochastic multiscale basis in Method S3
is reused multiple times. In particular, recalling the leading costs
\eqref{S2-cost} and \eqref{S3-cost}, we note that $\nreal = 2^{11}$
and $\nreal(2) = 2^9$, resulting in a global to local ratio
$\nreal/\nreal(2) = 4.0$. However, the reduction in runtime in Method
S3 compared to Method S2 is less significant, due to the additional
cost of inter-processor communication at each CG iteration. Similar
conclusions can be made in the sparse grid collocation case, where
Method S2 requires about 12\% of the subdomain solves and about 14\%
of the runtime of Method S1. Comparing Method S2 and Method S3, we
note that Method S3 requires about 52\% of the subdomain solves of
Method S2, which is consistent with the global to local ratio
$\nreal/\nreal(j) = 101/51 = 1.99$. We conclude that higher global to local
ratio results in increased reuse of the stochastic basis and increased
reduction in the number of subdomain solves.

\subsection{Test case 2}

In this example we test a porous medium with an inclusion of higher
and more heterogeneous permeability. The global domain $(0,1)^2$ is
divided in half with the Stokes region on the top and the Darcy region
on the bottom, and is partitioned into 32 subdomains. The outside boundary
conditions are as in Case 1, with the exception that the zero normal
velocity on the top edge is replaced with zero normal stress.

\begin{figure}
\centerline{\includegraphics[width=2.6in]{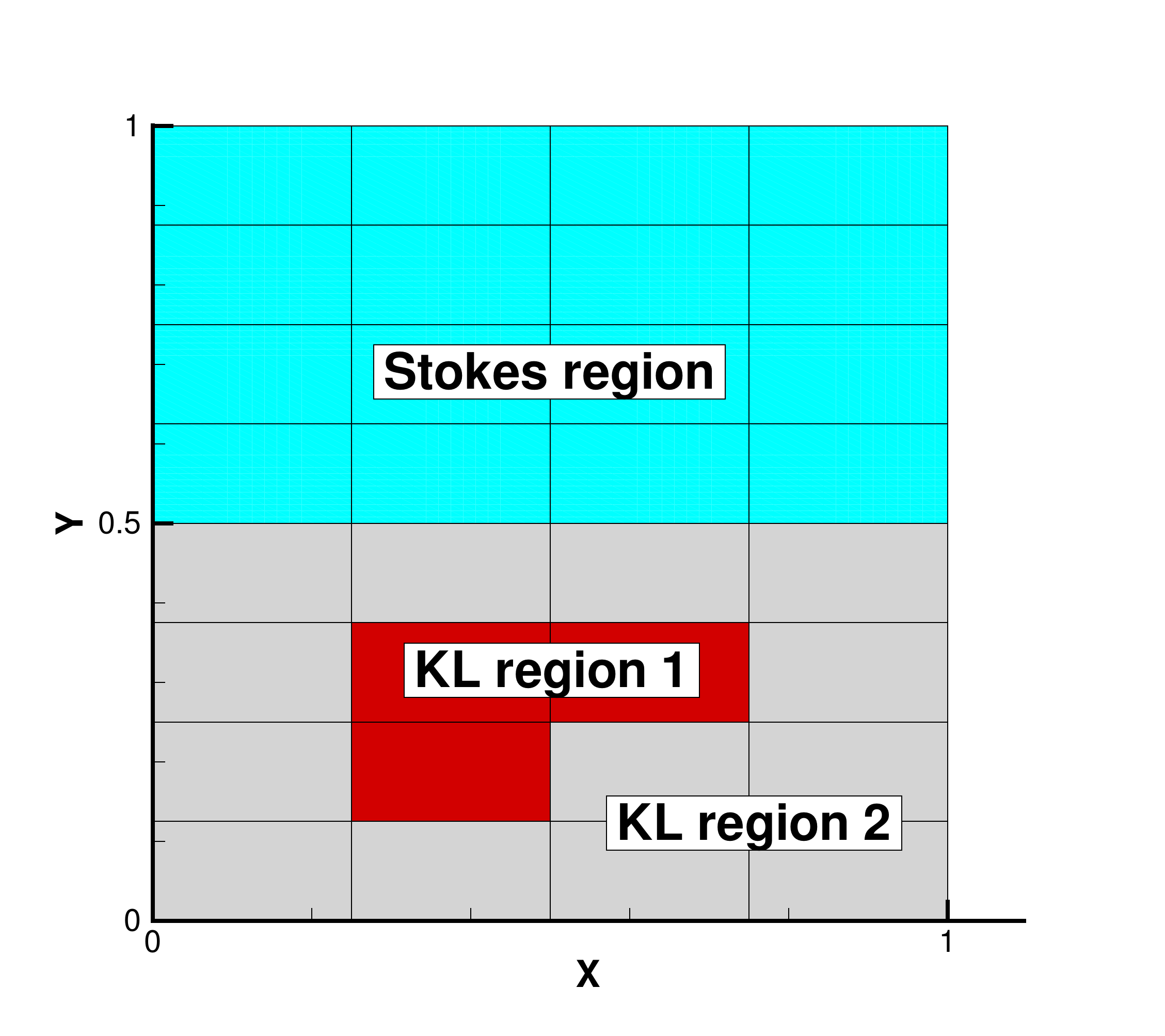}}
\caption{Case 2, illustration of KL regions.}
\label{Case2regions}
\end{figure}

\begin{figure}
	%
	\centerline{
		\includegraphics[width=.45\linewidth]{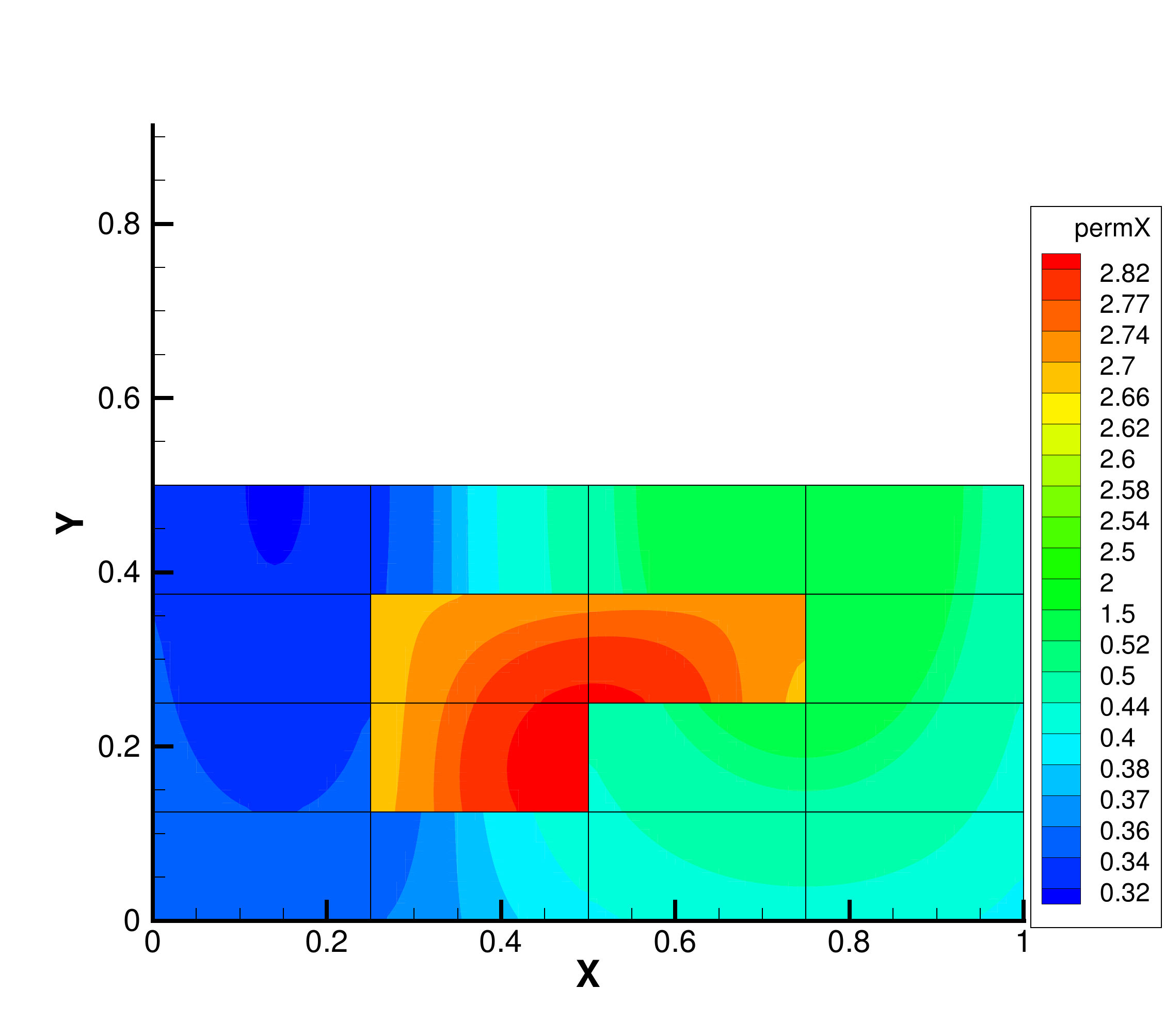}
		\includegraphics[width=.45\linewidth]{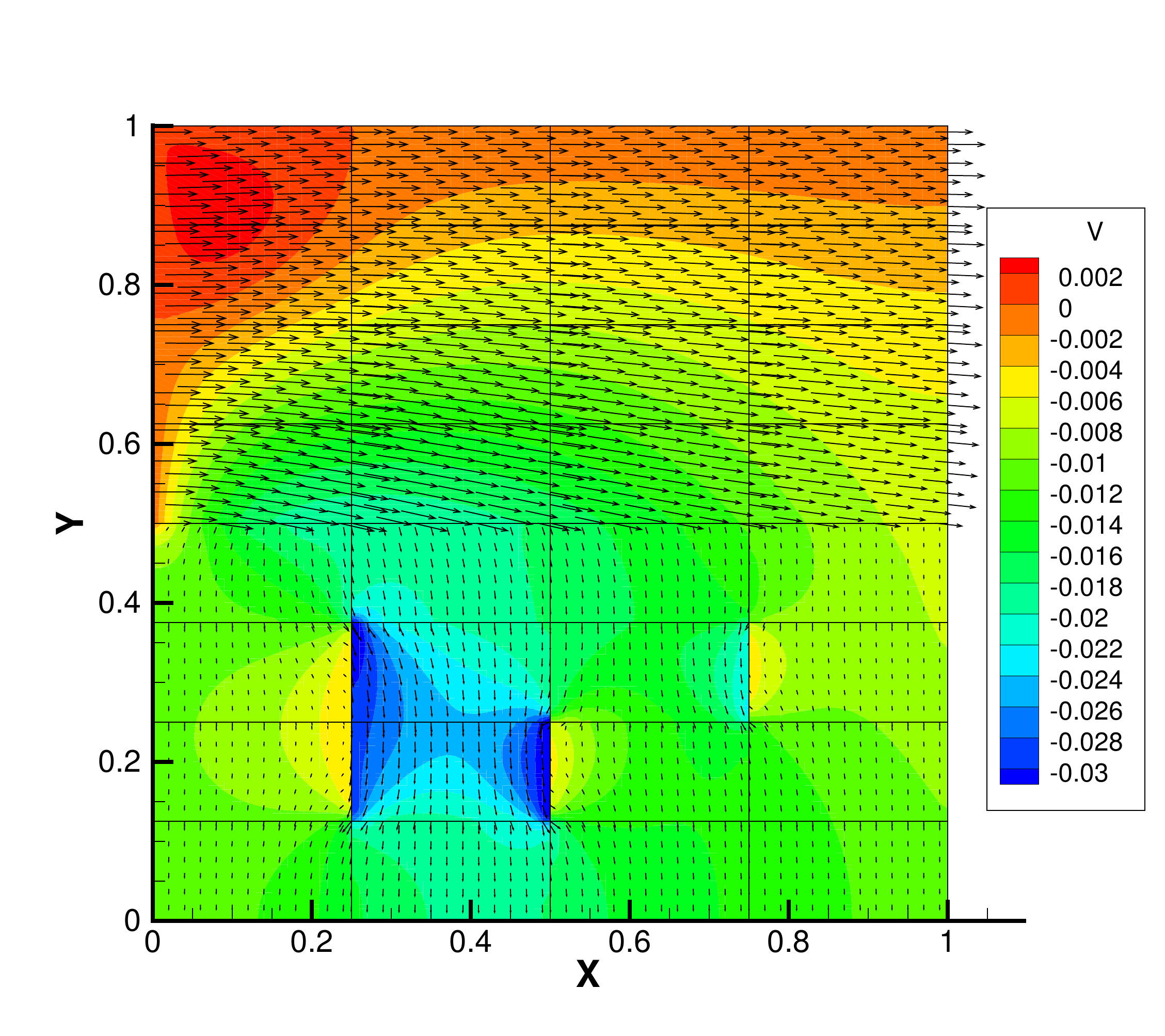}
	}
	\centerline{
		\includegraphics[width=.45\linewidth]{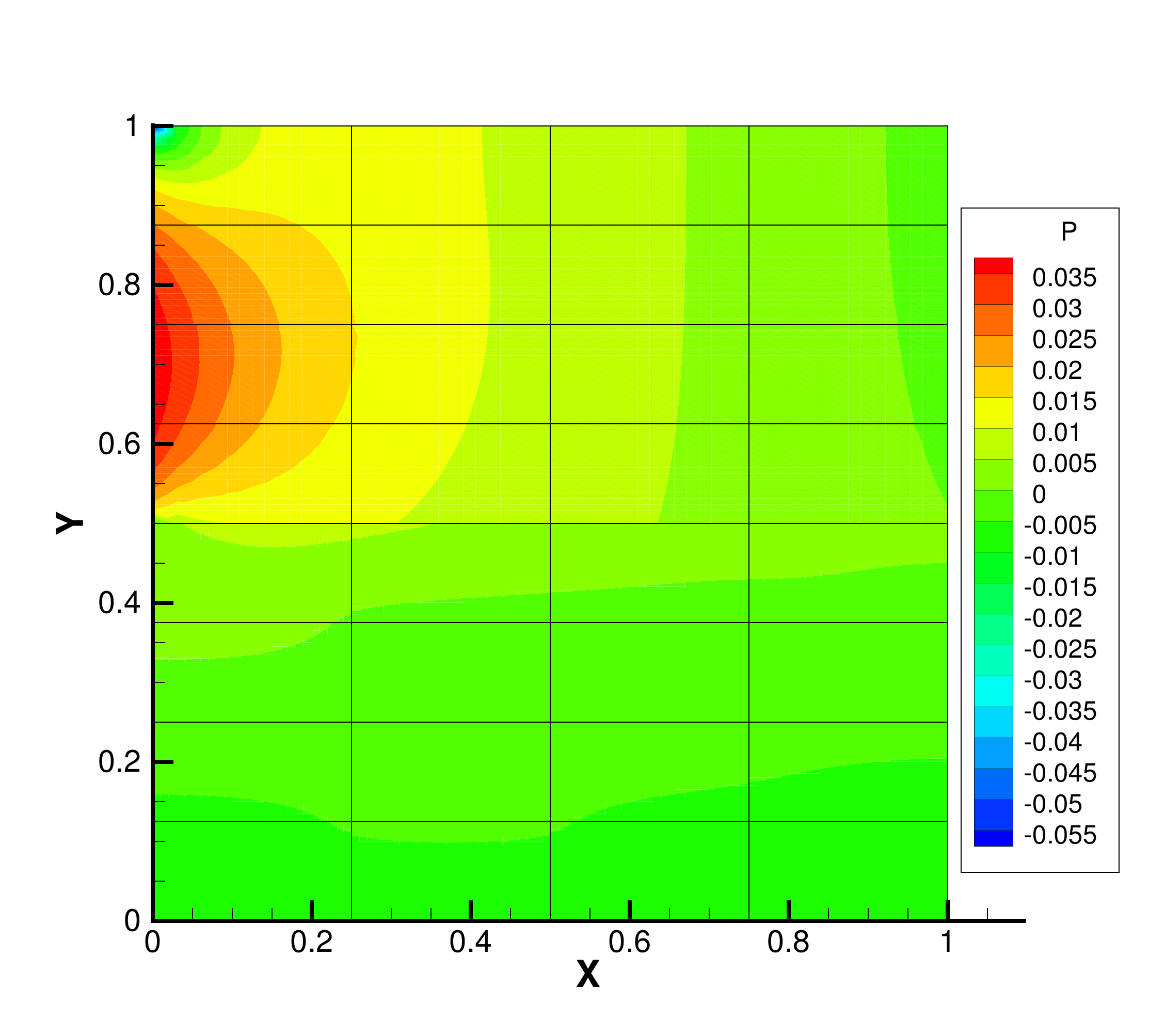}
		\includegraphics[width=.45\linewidth]{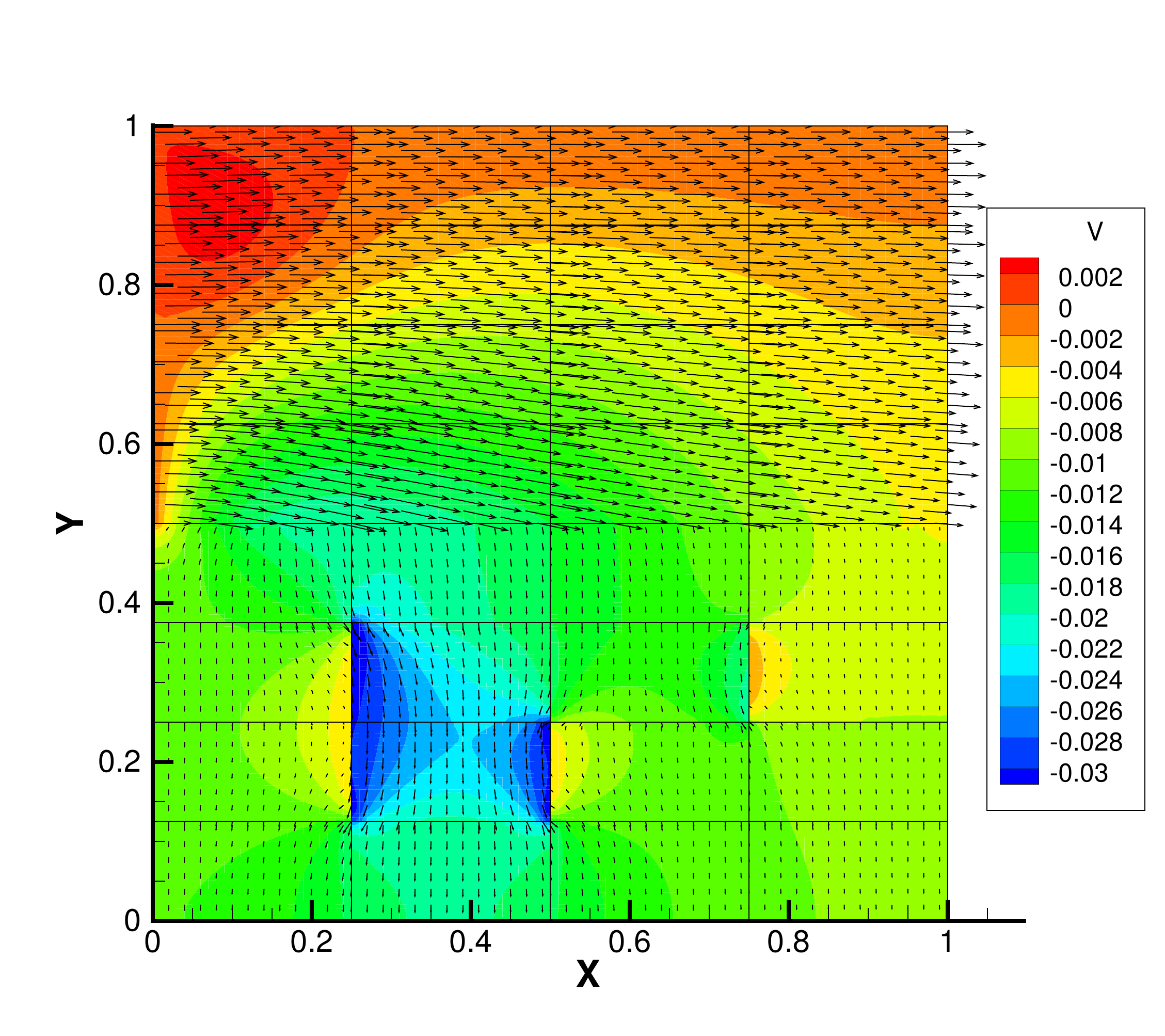}
	}
	\centerline{
		\includegraphics[width=.45\linewidth]{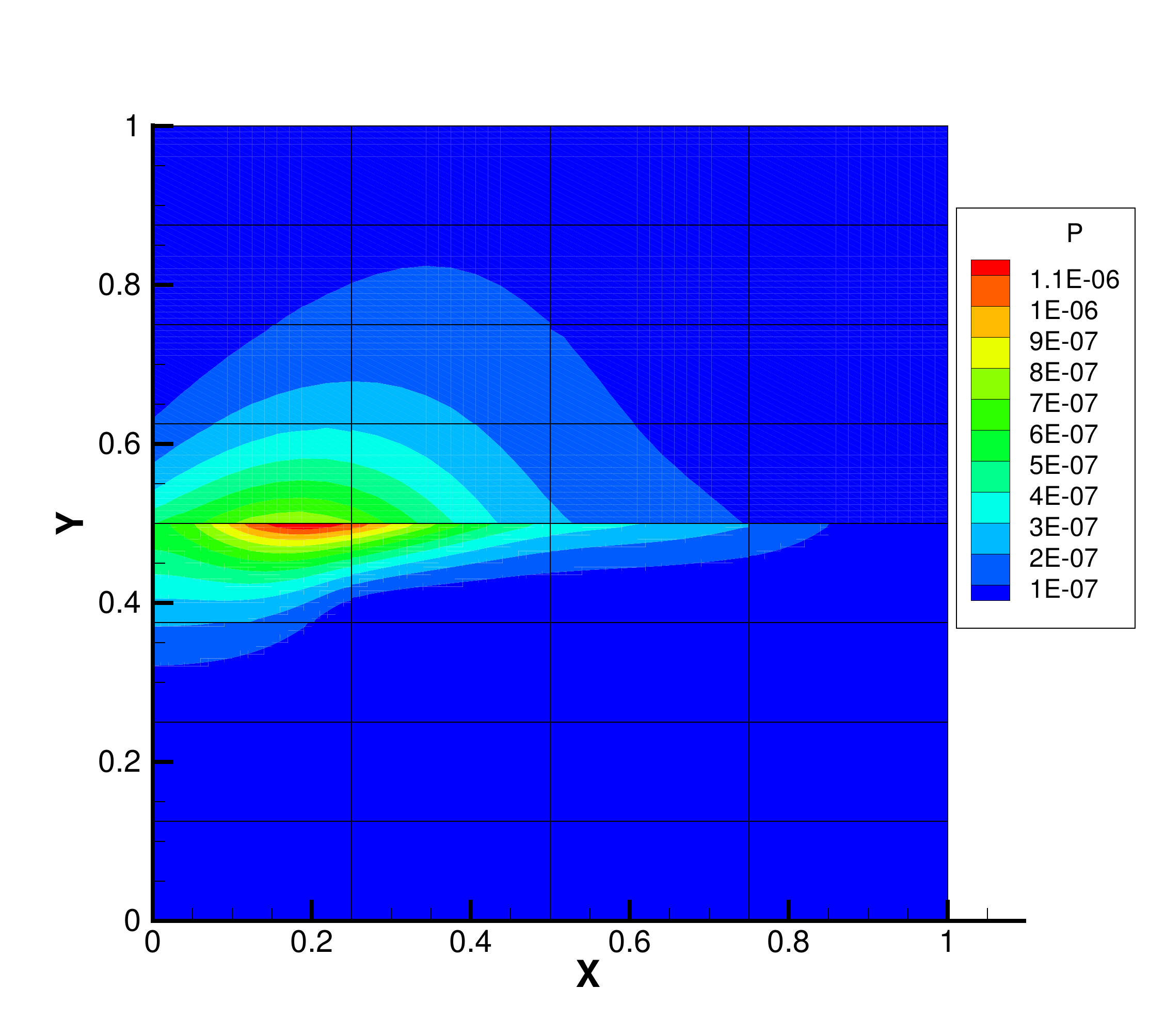}
		\includegraphics[width=.45\linewidth]{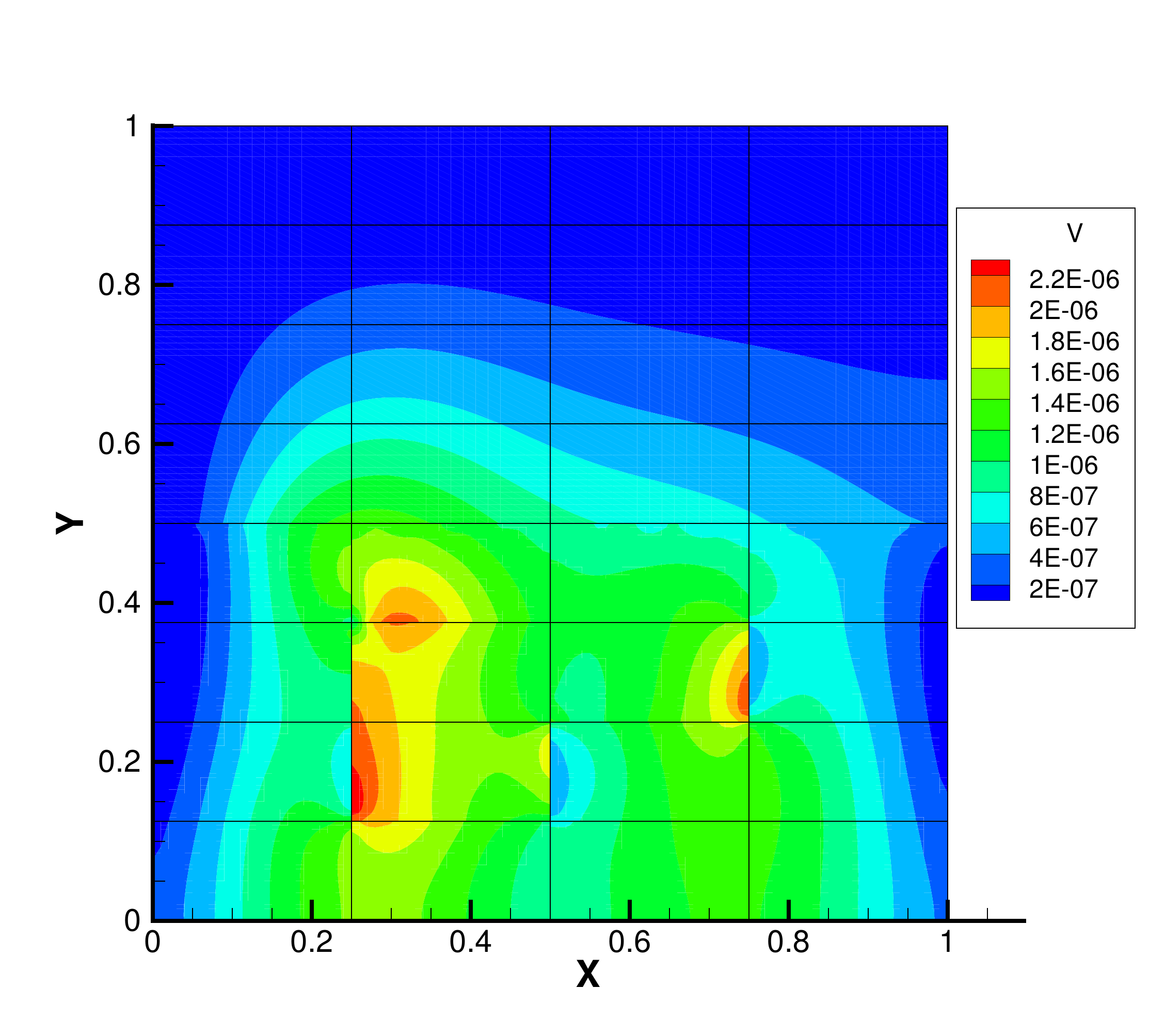}
	}
\caption{Case 2, tensor product grid, permeability realization 
  (top-left), velocity realization (top-right), pressure mean value
  (middle-left), velocity mean value 
  (middle-right), pressure variance (bottom-left), vertical velocity
  variance of (bottom-right). }
	\label{case2tensor}
\end{figure}

\begin{figure}
	%
	\centerline{
		\includegraphics[width=.45\linewidth]{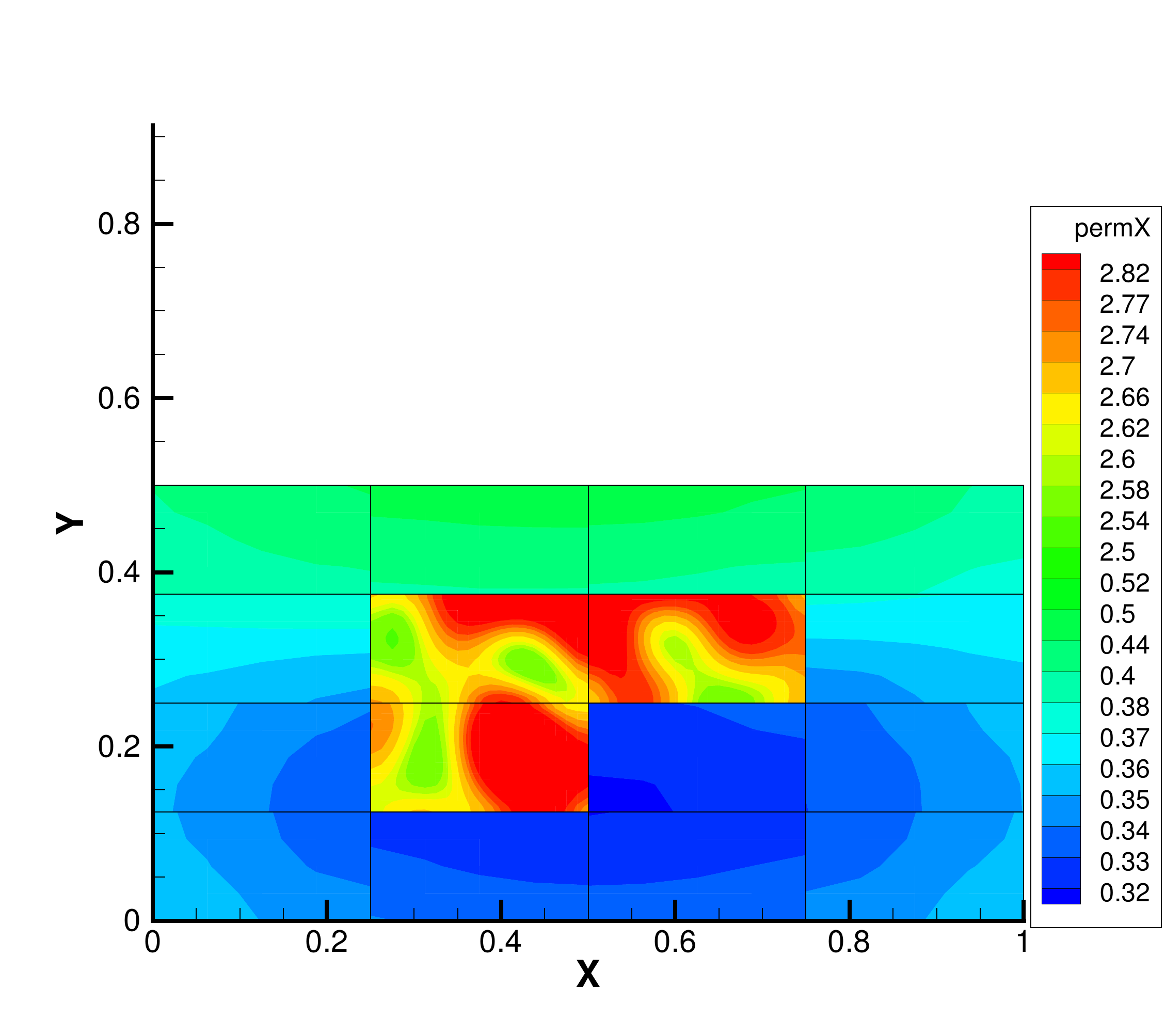}
		\includegraphics[width=.45\linewidth]{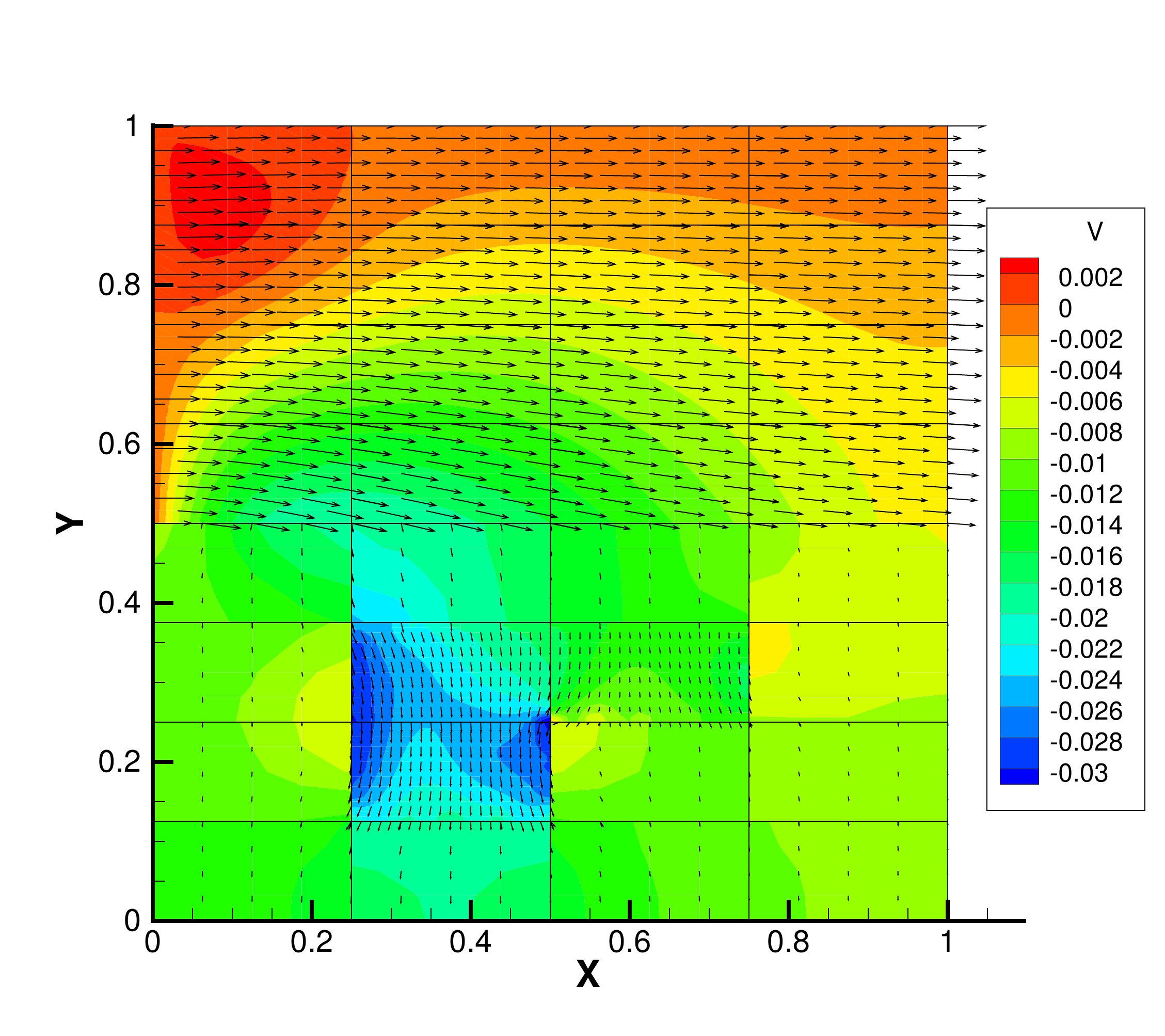}
	}
	\centerline{
		\includegraphics[width=.45\linewidth]{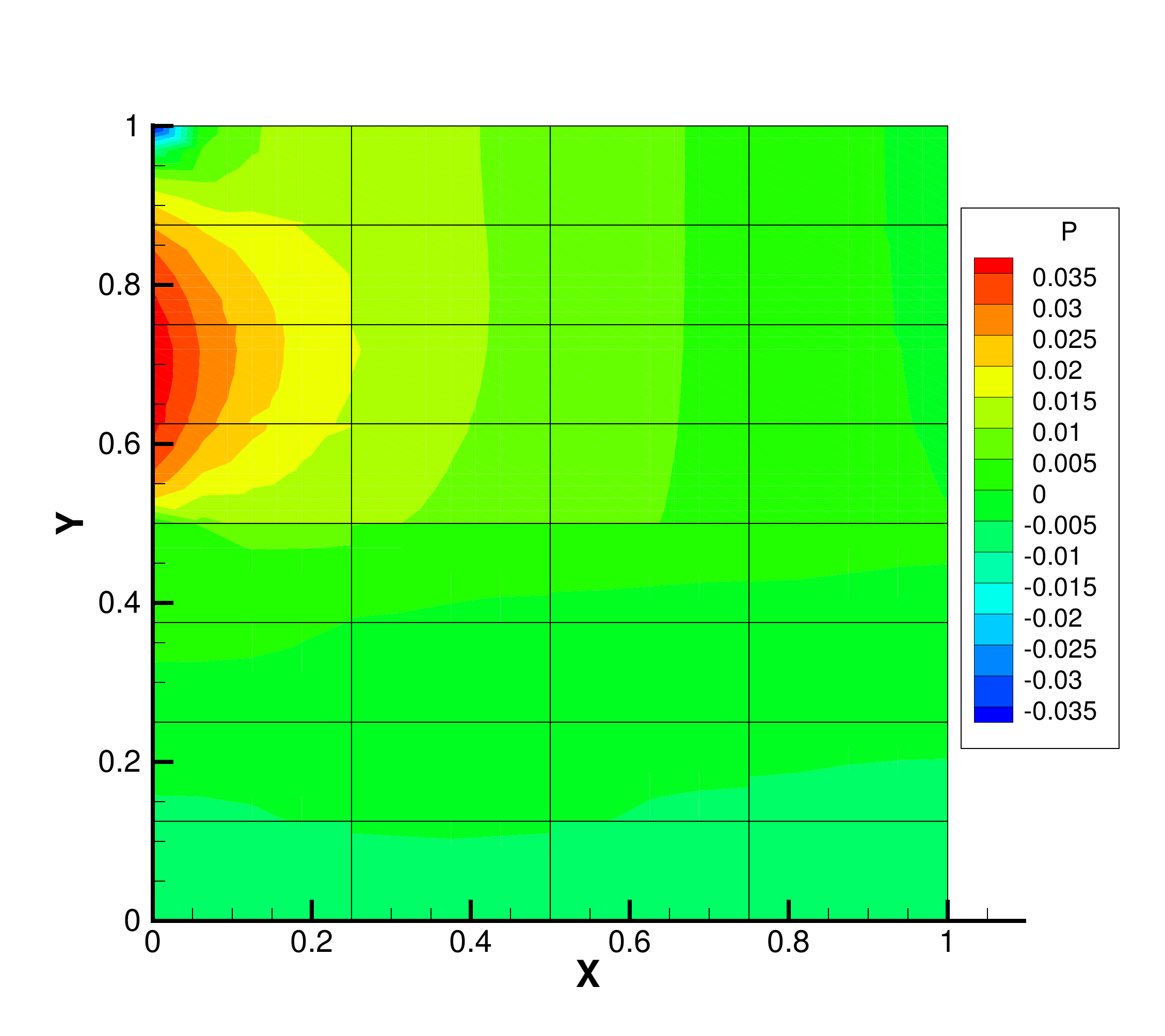}
		\includegraphics[width=.45\linewidth]{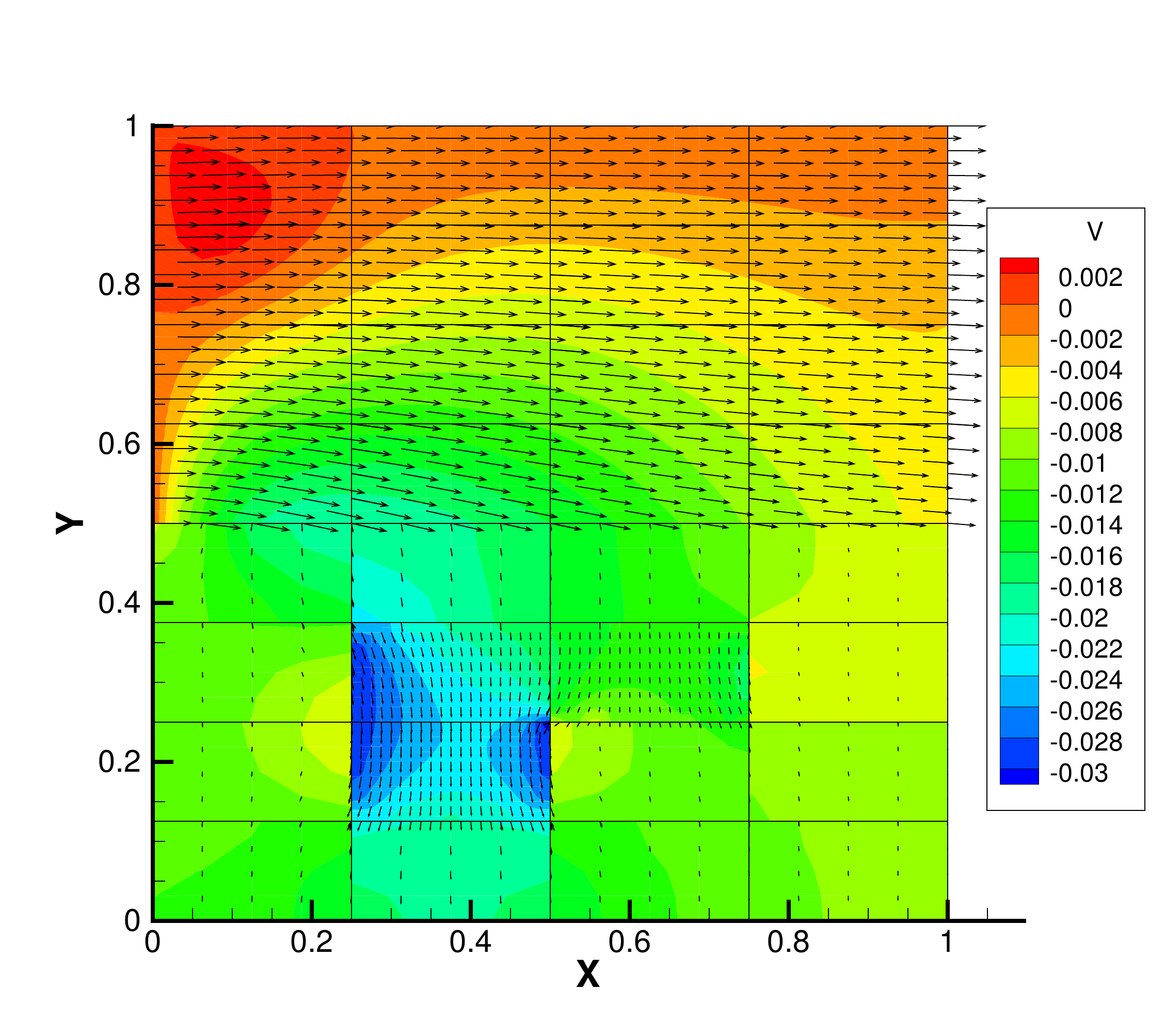}
	}
	\centerline{
		\includegraphics[width=.45\linewidth]{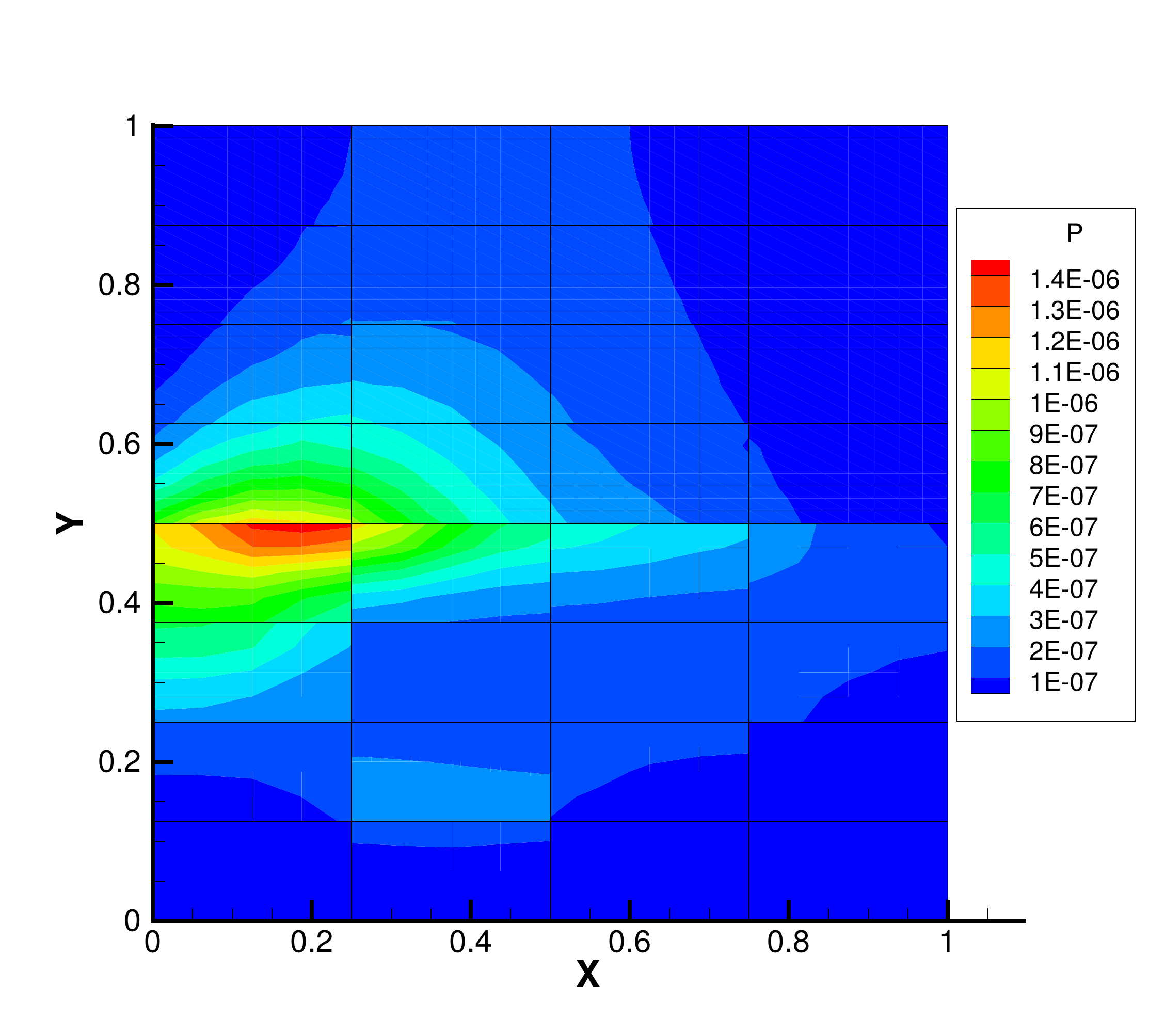}
		\includegraphics[width=.45\linewidth]{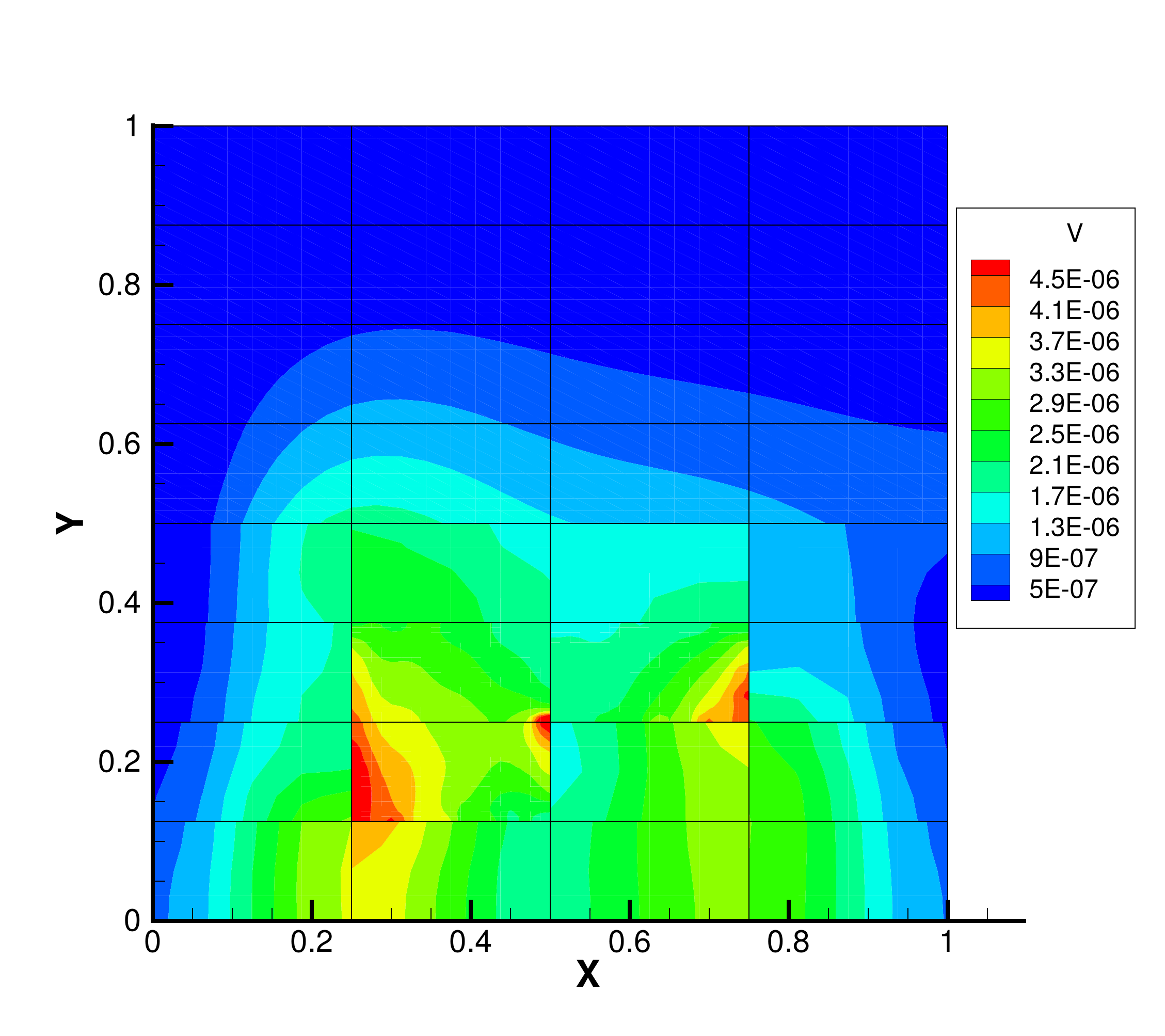}
	}
\caption{Case 2, sparse grid, permeability realization 
  (top-left), velocity realization (top-right), pressure mean value
  (middle-left), velocity mean value 
  (middle-right), pressure variance (bottom-left), vertical velocity
  variance of (bottom-right). }
\label{case2sparse}
\end{figure}

As shown in Figure \ref{Case2regions}, an L-shaped KL region 1 (red
part) is contained within the Darcy domain, with its complement
defined to be KL region 2 (gray part).  In KL region 1, the mean
permeability value is $e$, $\sigma_Y^{(1)} = 1.0$ and $\eta_j^{(1)}
= 0.01$.  In KL region 2, the mean permeability value is $e^{-1}$,
$\sigma_Y^{(2)} = 1.0$, and $\eta_j^{(2)} = 0.1$. Note that in the
inclusion, the mean permeability is more than 7 times larger and the
correlation lengths are 10 times smaller than in the rest of the Darcy
domain, see also the permeability realizations in top left panels in
Figures~\ref{case2tensor} and \ref{case2sparse}.  For the tensor
product collocation, the grid is isotropic with $\ncoll=2$ and
$\nterm=11$, where $\nterm(1)=3\times3$ and $\nterm(2)=2\times1$. For
the sparse grid, we have $\lmax=1$ and $\nterm=200$, where
$\nterm(1)=14\times14$ and $\nterm(2)=2\times2$.  We note that in both
collocation methods the number of KL terms is higher in the more
heterogeneous L-shaped inclusion. Furthermore, with 196 KL terms for
the sparse grid, the inclusion permeability is much more heterogeneous
compared to the tensor product grid. The number of global realizations is
$\nreal=2^{11}=2048$ for tensor product and $\nreal=2*200+1=401$
for sparse grid.

The physical space discretization is set as follows: for 
tensor product grid, we use
a local mesh of $25 \times 25$ in every Darcy subdomain and $8
\times 8$ in every Stokes subdomain. A mortar mesh of $10 \times 1$ is
used on every interface on $\Gm_{dd}$ and $4\times 1$ on every interface on 
$\Gm_{sd}\cup\Gm_{ss}$; for sparse grid, we use a local mesh of $20
\times 20$ in every subdomain in the inclusion and $4 \times 4$ elsewhere. A
mortar mesh of $10 \times 1$ is used on every interior interface of
the inclusion, $4\times 1$ on each of its outside boundaries,
and $2\times 1$ elsewhere.

Computed realizations, mean, and variance of the solution are shown in
Figures \ref{case2tensor} and \ref{case2sparse} for tensor product and
sparse grid, respectively. The general behavior of the flow is similar
to Case 1, with a channel-like flow in the Stokes region infiltrating
downward into the Darcy region. The largest vertical velocity is in
the high permeability L-shaped inclusion. The computed solution
variance indicates that the largest uncertainty in the velocity is
within the inclusion, while the largest pressure uncertainty is along
the Stokes-Darcy interface. As expected, the solution uncertainty is
higher in the Darcy region, since there is no stochastic parameter in
the Stokes region. Nevertheless, the solution uncertainty extends well
into the Stokes region, being highest along the Stokes-Darcy
interface. This is due to the fact that the coupling conditions imply
that the solution is stochastic in the entire Stokes-Darcy domain.

The maximum number of subdomain solves and runtimes for Case 2 are
shown in Table~{\ref{Case2}}. We note that in the multiscale flux
basis Methods S2 and S3, the subdomain with the maximum number of
solves is the one at the corner of L-shaped region, since it has the
highest number of mortar degrees of freedom on its interfaces. For the
tensor product grid, Method S2 requires only 11\% of the subdomain
solves in Method S1 and 52\% of the runtime. These numbers for the
sparse grid are 13\% and 83\%, respectively. Again, the smaller gain
in runtime is due to the inter-processor communication overhead, which
is even more significant in the example because of the larger number
of subdomains. As in Case 1, the stochastic multiscale basis Method S3
results in additional computational gain compared to Method S2. For
tensor product grid Method S3 requires 26\% of the subdomain solves of
Method S2, which is consistent with the global to local ratio
$\nreal/\nreal(1) = 2^{11}/2^{9} = 4.0$. However, the gain with sparse
grid is much smaller, with Method S3 requiring 98\% of the maximum
number of subdomain solves of Method S2. This again can be explained
with the global to local ratio. In particular, $\nreal/\nreal(1) =
401/193 = 1.02$, thus the stochastic basis in KL region 1 is reused
very few times. Such situation occurs when most of the dimensions in
the global stochastic space are contained within one KL region, as is
the case here, with 196 KL terms in KL region 1 out of total of
200 KL terms. 

\begin{table}
\begin{center}
{\em \textbf{ Tensor Product} Collocation, $\nterm=11$, $\ncoll=2$ (2048 realizations)
}
\begin{tabular}{||l||c|c|c||}
\hline
& \textbf{Method S1} & \textbf{Method S2} & \textbf{Method S3} \\
\hline
Max. number of solves  & 3,033,494 & 331,776 & 86,016 \\ 
\hline
Runtime in seconds   & 5893.3   & 3155.9   & 2707.6  \\
\hline
\end{tabular}
		
\vspace{1.5ex}
{\em \textbf{Sparse Grid} Collocation,  $\nterm=200$, $\lmax=1$ (401 realizations) 
}
\begin{tabular}{||l||c|c|c||}
\hline
& \textbf{Method S1} & \textbf{Method S2} & \textbf{Method S3}  \\
\hline
Max. number of solves  & 346,709 & 45,714 & 44,818  \\
\hline
Runtime in seconds   & 525.5   & 435.1   & 391.2    \\ 
\hline
\end{tabular}
\end{center}
\caption{Number of subdomain solves and runtime in seconds for Case 2.}
\label{Case2}
\end{table}

\subsection{Test case 3}

\begin{figure}
\centering
\begin{minipage}{.46\textwidth}
  \noindent
    \includegraphics[width=2.7in]{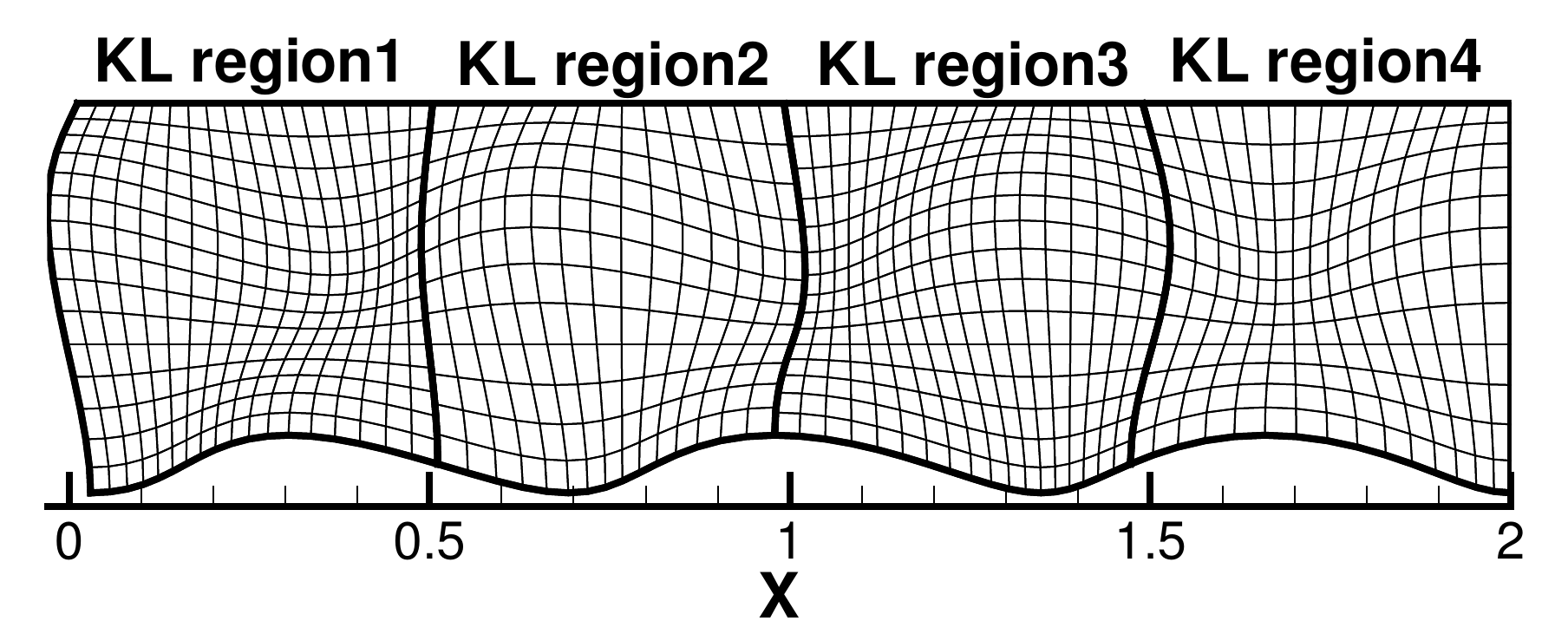}\\
      \includegraphics[width=2.95in]{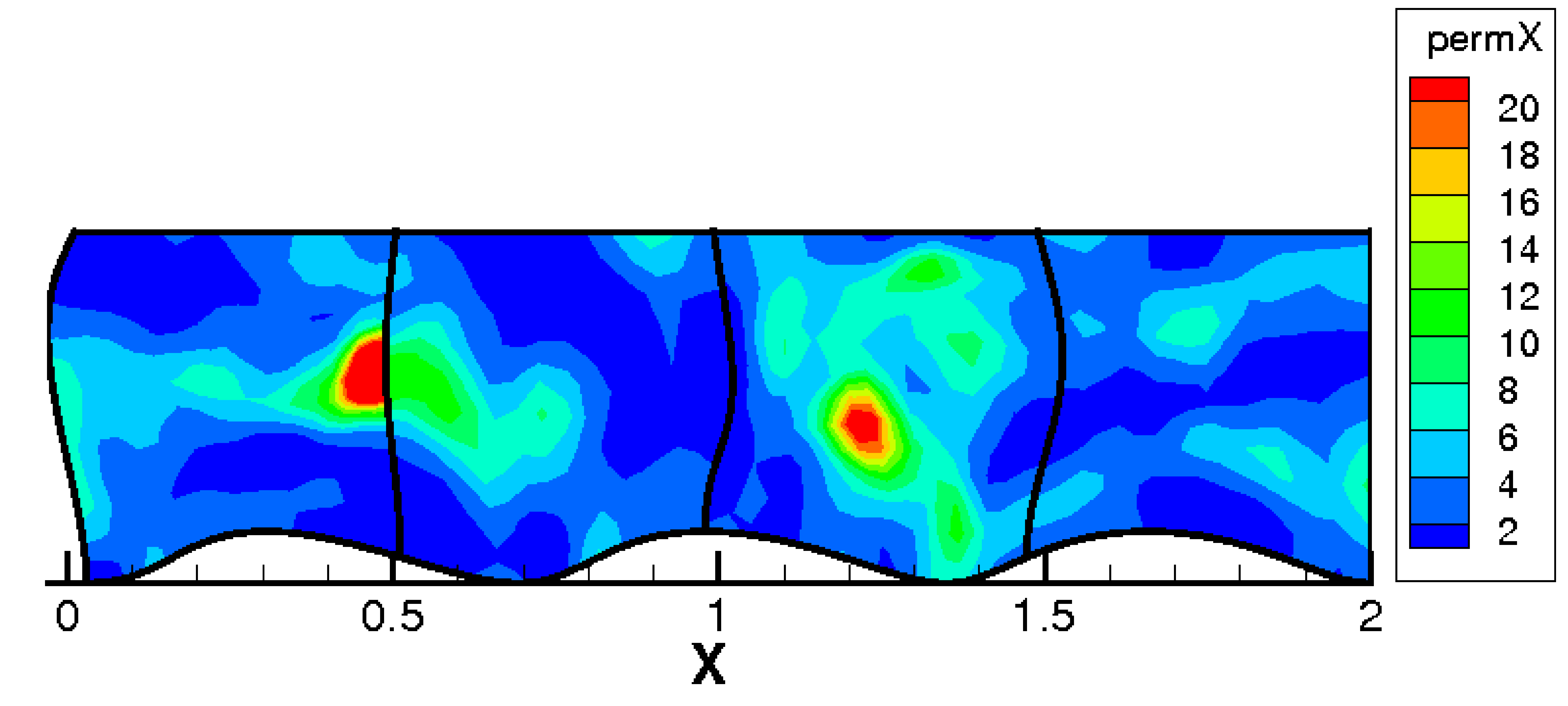}
\end{minipage}
\hfill
\begin{minipage}{.5\textwidth}
    \includegraphics[width=.99\textwidth]{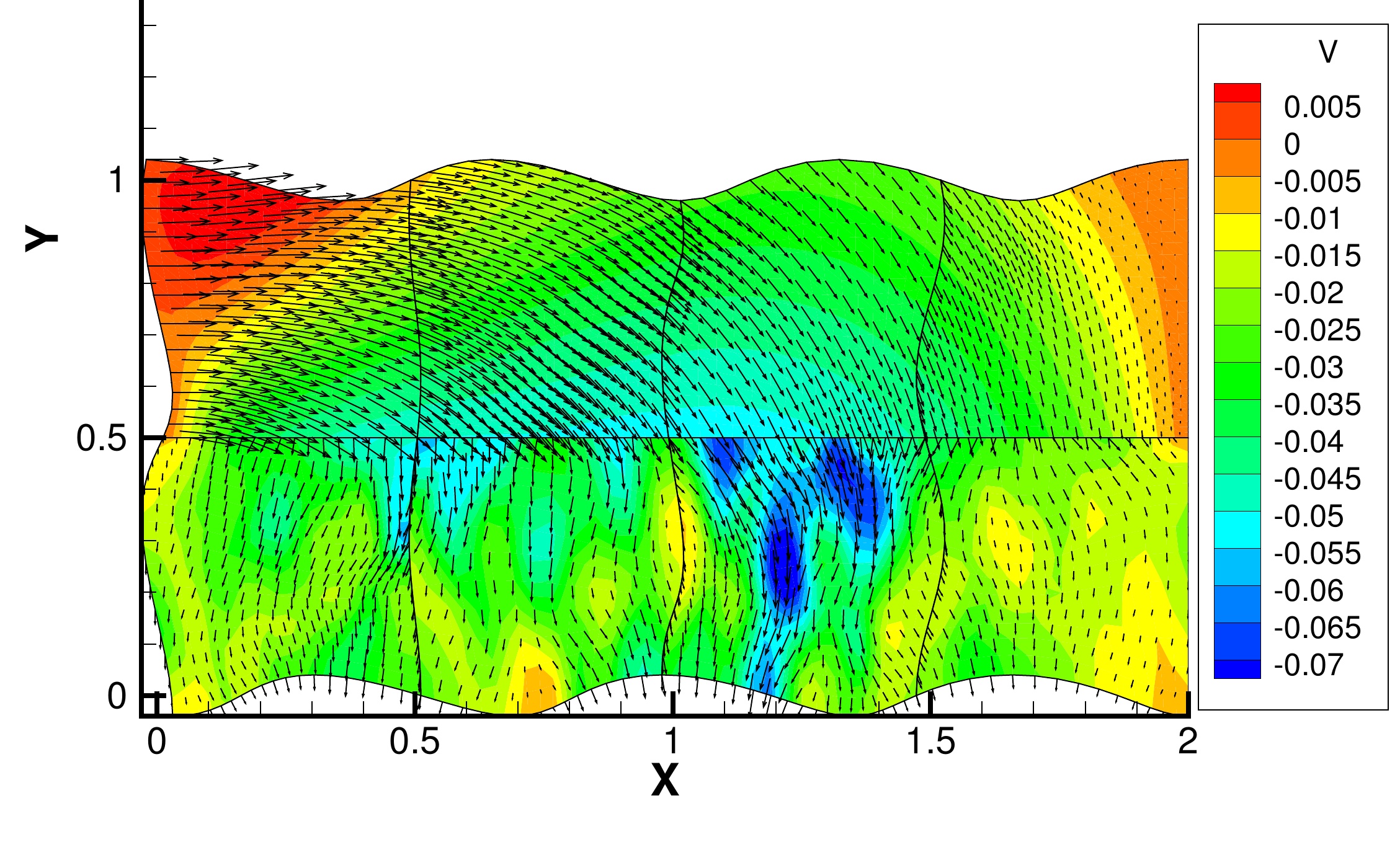}
\end{minipage}
\caption{Case 3, KL regions and meshes (top-left), permeability
  realization (bottom-left), velocity realization (right).}
\label{Case3region}
\end{figure}

The third test is a coupled surface water and groundwater flow with
realistic geometry and highly heterogeneous permeability. The
computational domain is shown in Figure~\ref{Case3region} (right).
There are $4\times 2=8$ subdomains on the global irregular region with
Stokes in the top half and Darcy on the bottom.  The outside boundary
conditions are given as follows: in Darcy, zero pressure is specified
on the bottom, and no-flow on left and right; for the Stokes region,
inflow velocity is specified on the left and zero velocity on the
right, along with zero normal and tangential stresses on the top. The
mean permeability over the entire Darcy region is generated by a
single realization of a global KL expansion truncated after 400 terms,
with mean value one and covariance given in \eqref{covariance} with
variance $\sigma_Y = 2.1$ and correlation lengths $\eta_1 = 0.1$ and
$\eta_2 = 0.05$.  There are four KL regions in Darcy, represented by
each one of the four Darcy subdomains, as shown in Figure
\ref{Case3region} (top-left).  In all KL regions, the variance is
$\sigma_Y^{(i)} = 1.0$ and the correlation lengths are $\eta_j^{(i)} =
0.1$. One realization of the permeability field with tensor product
collocation is presented in Figure \ref{Case3region} (bottom-left).
For tensor product collocation, the grid is isotropic with $\ncoll=2$,
the stochastic dimensions are set by $\nterm(j)=2\times 1$ for
$j=1,2,3,4$ and $\nterm=8$ in total. For sparse grid collocation, we
have $\lmax=1$ and $\nterm=100$ with $\nterm(j)=5\times 5$ for
$j=1,2,3,4$. The physical grids are alternating $18\times 15$ and
$15\times 12$ in the Darcy subdomains, and $12\times 15$ and $9\times
12$ in the Stokes subdomains.  Mortar meshes of $4\times1$ are used on
all interfaces.

To handle the irregular geometry, we employ the multipoint flux mixed
finite element method on quadrilaterals in the Darcy region
\cite{WXY-12}, which can be reduced to a cell-centered finite
difference scheme for the pressure.  We impose the mortar conditions
on curved interfaces by mapping the physical grids to reference grids
with flat interfaces.  For more details about this implementation, the
reader is referred to \cite{Song-Wang-Yotov}.

\begin{figure}
  \centerline{
  \hspace{.3cm}   
  \includegraphics[width=.44\linewidth]{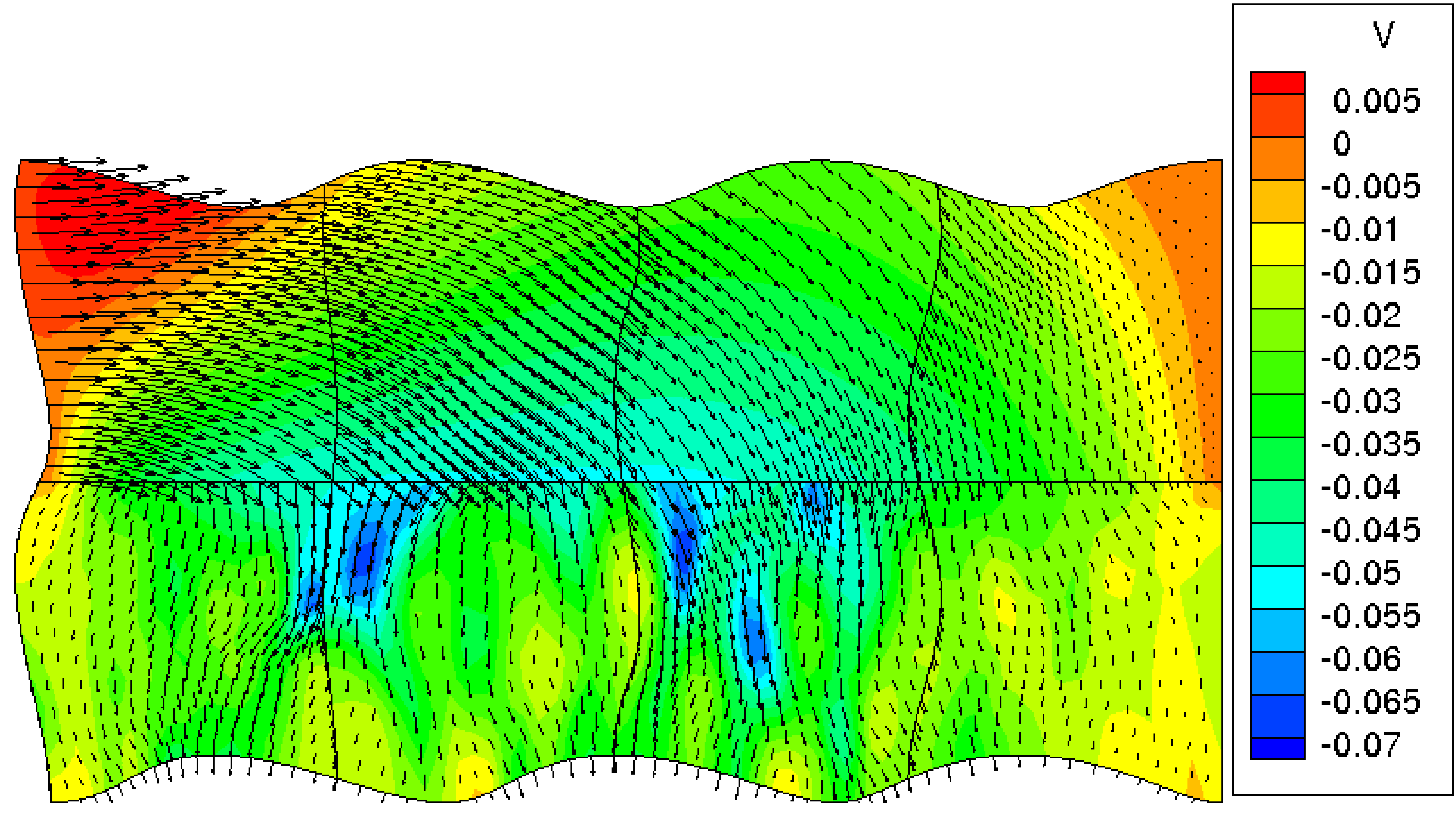}
\hspace{.9cm}  
\includegraphics[width=.44\linewidth]{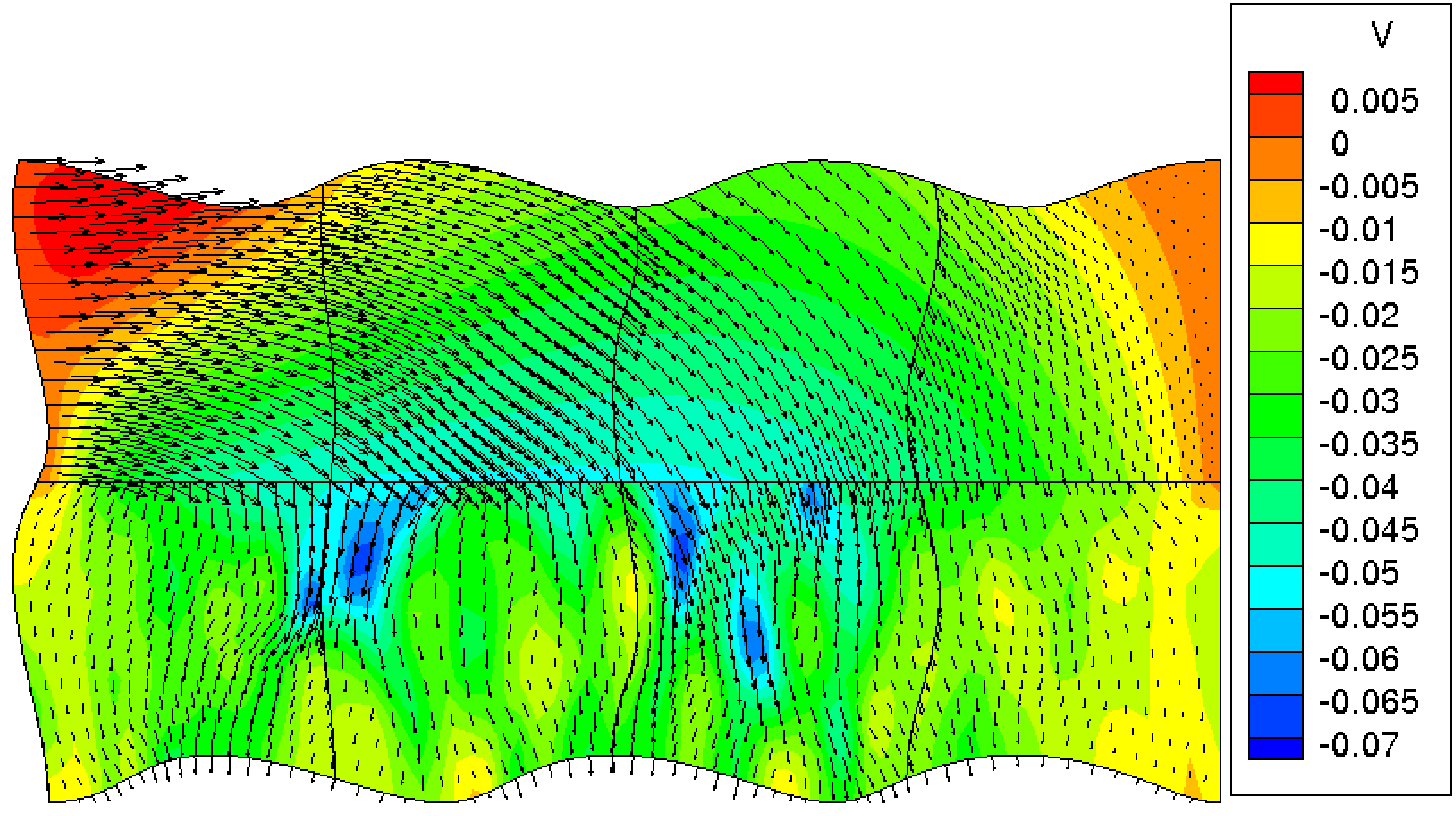}
}
\centerline{
\includegraphics[width=.5\linewidth]{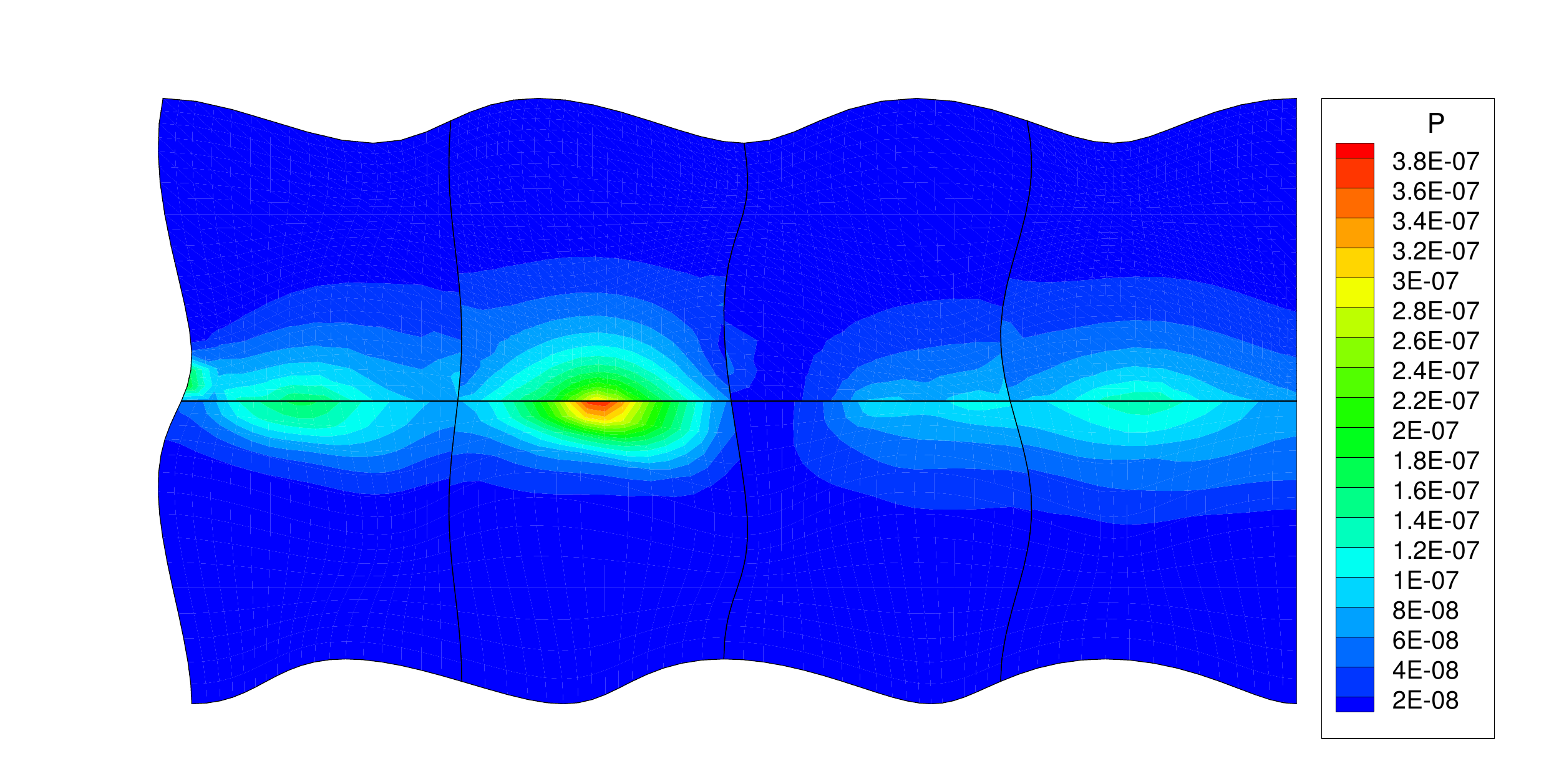}
\includegraphics[width=.5\linewidth]{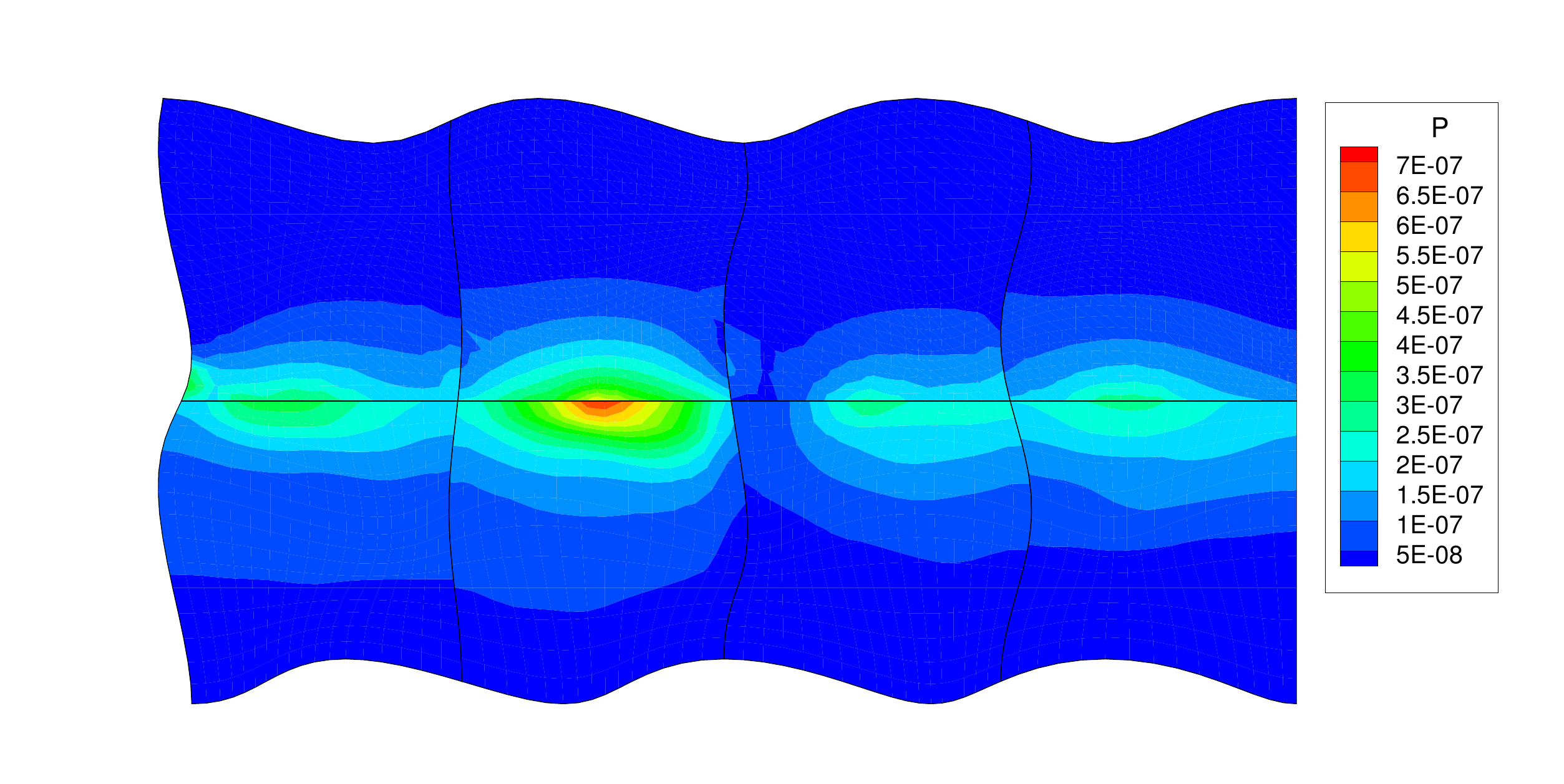}
}
\centerline{
\includegraphics[width=.5\linewidth]{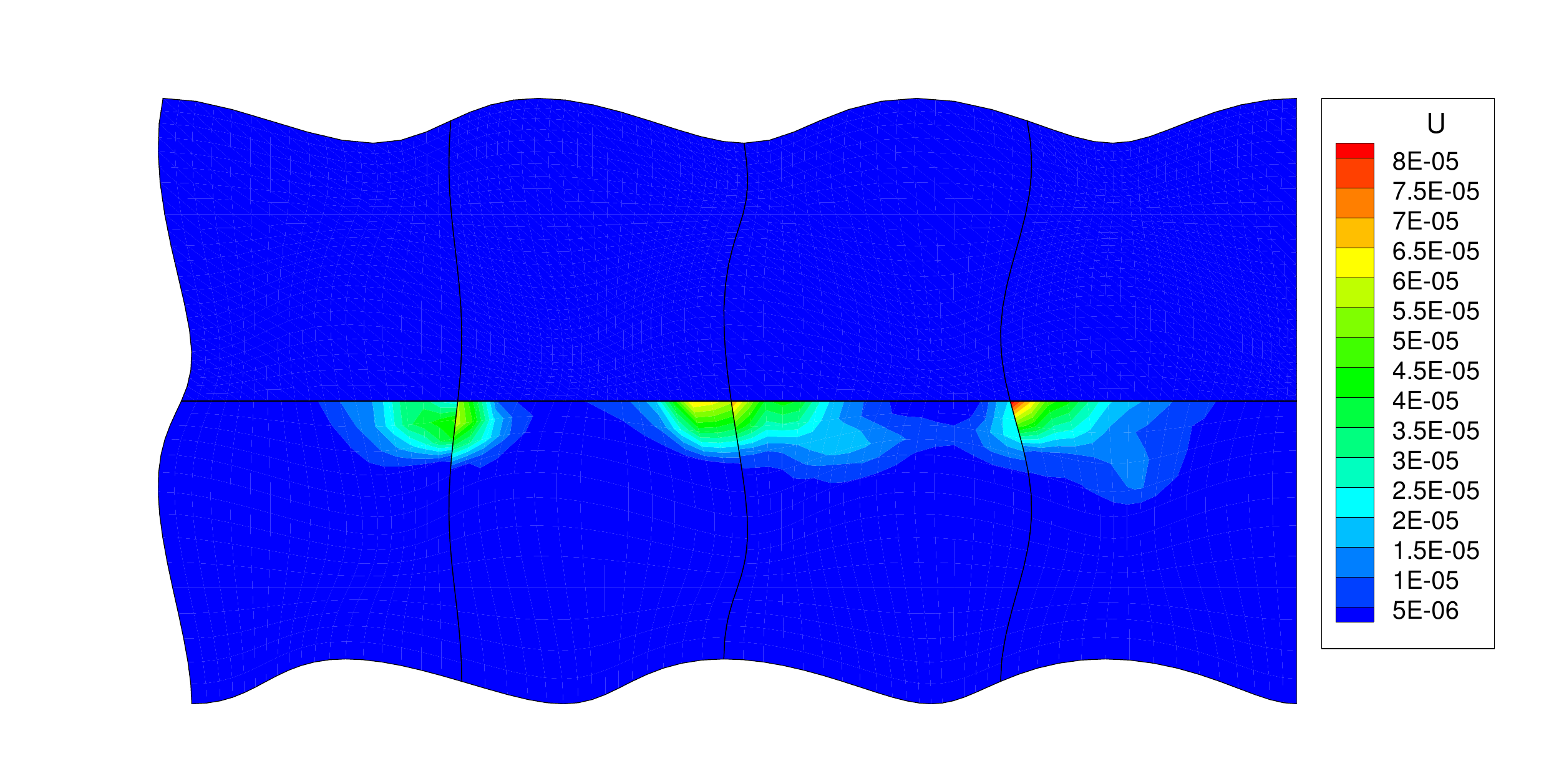}
\includegraphics[width=.5\linewidth]{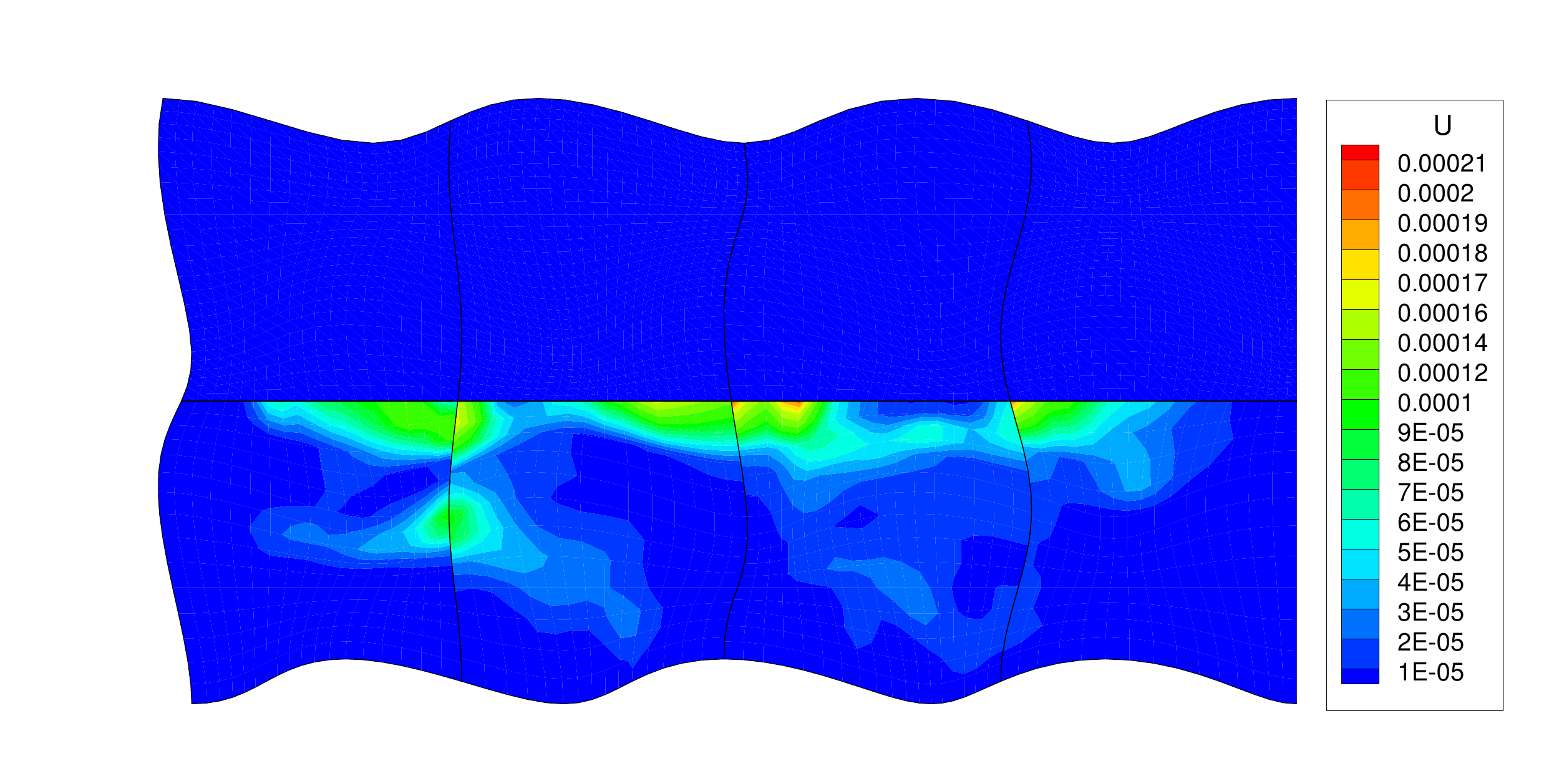}
}
\centerline{
\includegraphics[width=.5\linewidth]{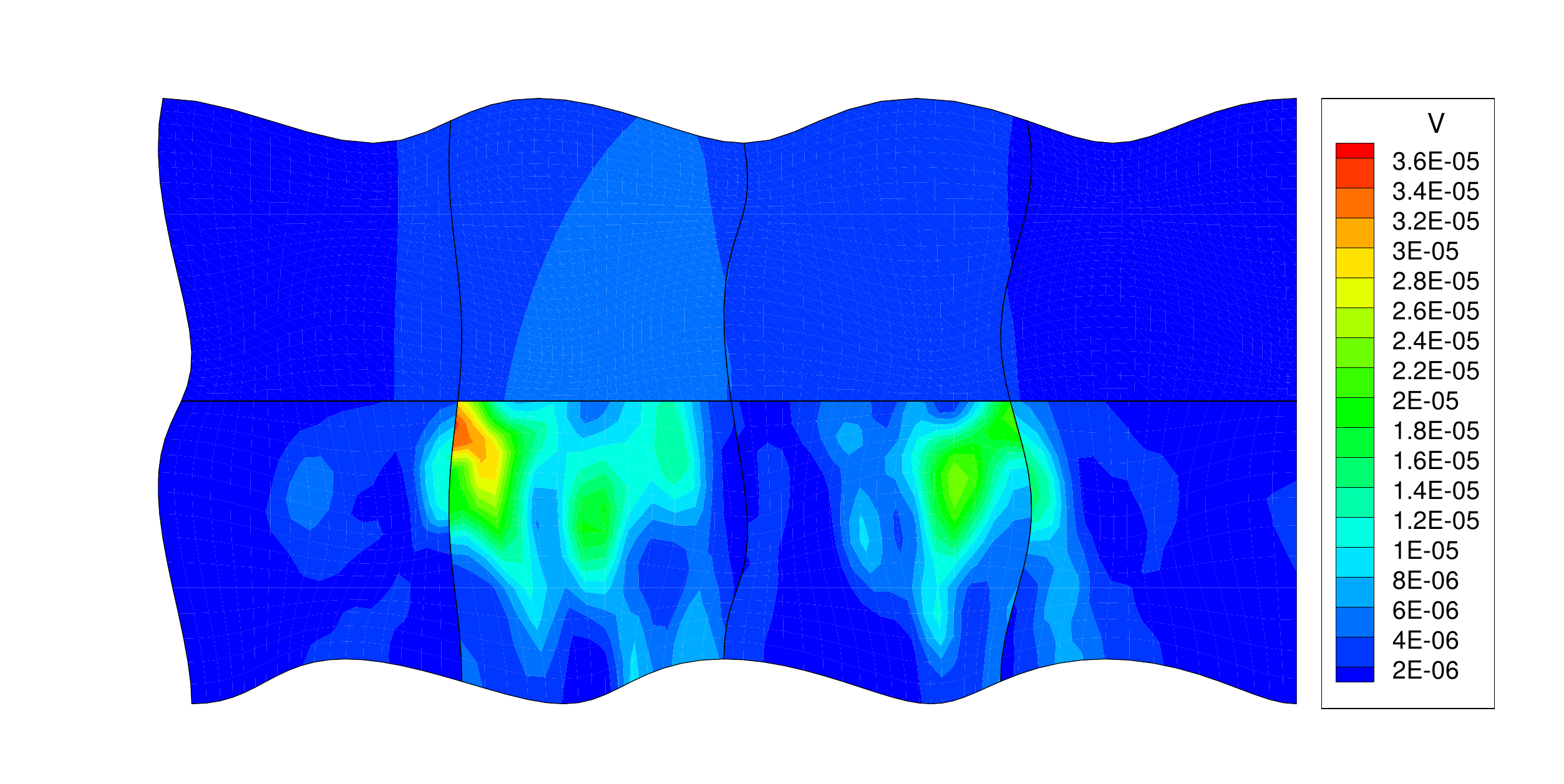}
\includegraphics[width=.5\linewidth]{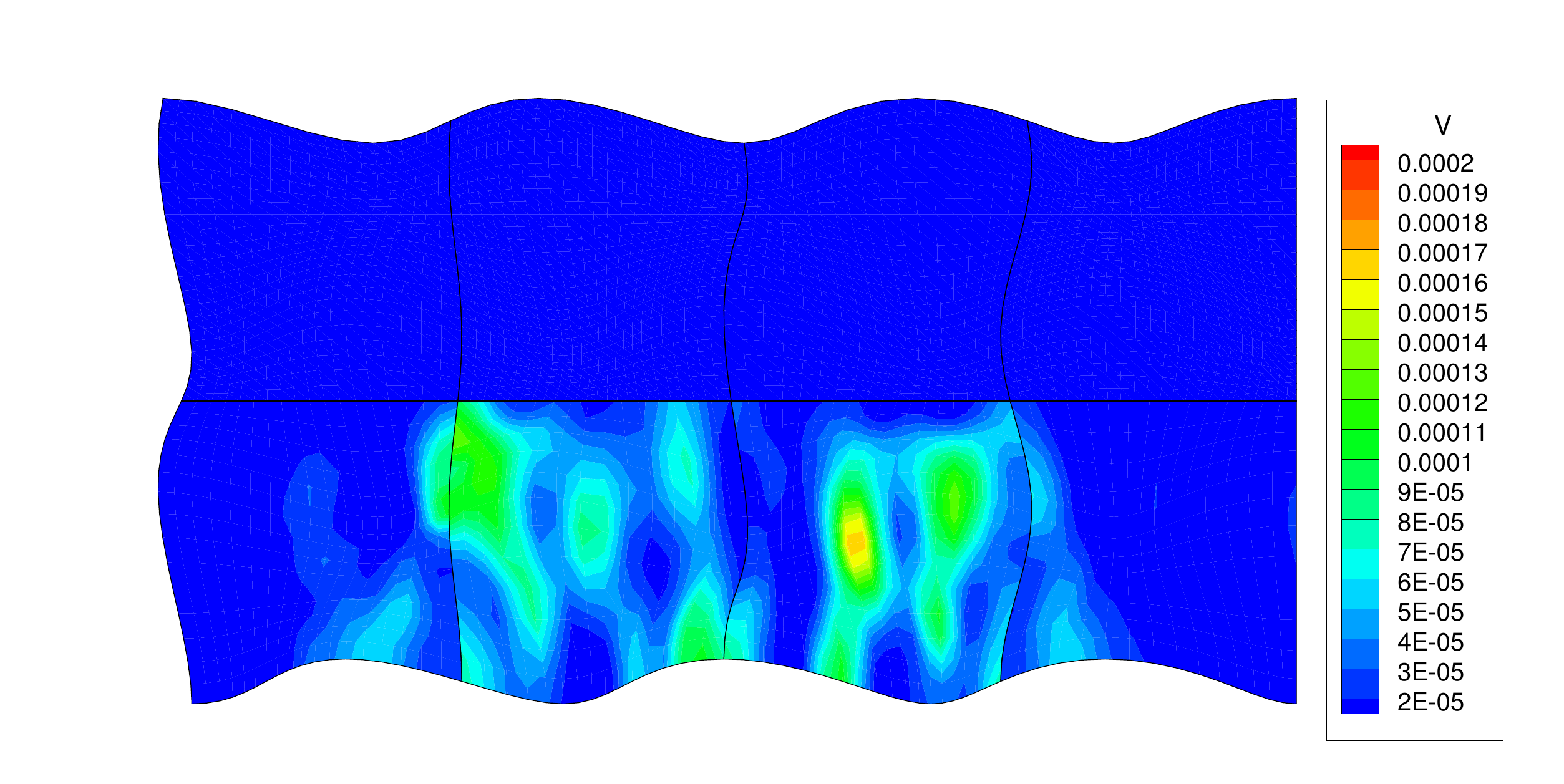}
}
\caption{Case 3, solution with tensor product grid (left) and sparse grid
  (right); from top to bottom: velocity mean, pressure variance,
  horizontal velocity variance, vertical velocity variance.}
\label{Case3-mean-var}
\end{figure}

Figure~\ref{Case3-mean-var} shows the mean and variance of the
computed solution with tensor product grid (left) and sparse grid
(right).  The predicted mean velocity with the two methods is similar,
indicating that it is well resolved with both stochastic collocation
methods.  The computed variance is larger and better resolved with the
sparse grid, due to the larger number of KL terms.

\begin{table}
\begin{center}

{\em  \textbf{ Tensor Product} Collocation, $\nterm=8$, $\ncoll=2$ (256 realizations)
}
\begin{tabular}{||l||c|c|c||}
\hline
& \textbf{Method S1} & \textbf{Method S2} & \textbf{Method S3} \\
\hline
Max. number of solves  & 58,728 & 21,248 & 960 \\ 
\hline
Runtime in seconds   & 255.5   & 192.8   & 91.7 \\
\hline
\end{tabular}

\vspace{1.5ex}
 {\em \textbf{Sparse Grid} Collocation,  $\nterm=100$, $\lmax=1$ (201 realizations) 
 }
\begin{tabular}{||l||c|c|c||}
\hline
& \textbf{Method S1} & \textbf{Method S2} & \textbf{Method S3} \\
\hline
Max. number of solves  & 45,031 & 16,683 & 3,051 \\
\hline
Runtime in seconds   & 178.3   & 107.4   & 82.6 \\ 
\hline
\end{tabular}
\end{center}
\caption{Number of subdomain solves and runtime in seconds for Case 3.}
\label{Case3}
\end{table}

The computational cost for Case 3 is presented in Table~\ref{Case3}.
The number of global stochastic realizations is $\nreal=256$ in tensor
product and $\nreal=201$ in sparse grid. The number of subdomain
solves for Method~S1 is smaller compared to the previous examples, due
to the fewer realizations and fewer CG iterations. Thus smaller
computational savings with Methods S2 and S3 are
expected. Nevertheless, with tensor product grid, Method~S2 requires
36\% of the subdomain solves of Method S1 and 75\% of the
runtime. With sparse grid, it requires 37\% of the solves of Method~S1
and 60\% of the runtime. The global to local ratio in this case is
larger, $\nreal/\nreal(j)=256/4=64$ for tensor product grid and
$\nreal/\nreal(j)=201:51= 3.94$ for sparse grid. This results in
significant savings in Method~S3 compared with Method~S2 in terms of
subdomain solves, with Method~S3 requiring only 4.5\% of the solves
and 48\% of the runtime with tensor product grid, and 18\% of the
solves and 77\% of the runtime with sparse grid.

\section{Conclusions}
We have presented three algorithms for solving the coupled stochastic
Stokes-Darcy problem with non-stationary porous media, one of which
does not use a multiscale flux basis and the other two use
deterministic and stochastic multiscale flux basis, respectively. The
discretization includes stochastic collocation methods in the
stochastic space and the MMMFEM in the physical space. The numerical
tests in Section \ref{numerical} demonstrate that the proposed
algorithms can compute accurate statistical moments of the solution
for a wide range of permeability heterogeneity and uncertainty, number
of KL regions and KL terms.  We generally observe a significant
reduction in the number of subdomain solves when the multiscale flux
basis is utilized for solving the coarse-scale mortar interface
problem. The stochastic multiscale flux basis results in significant
additional savings, due to its reuse during the global collocation
loop, especially in the cases with large global to local realizations
ratio, such as using tensor product grid collocation or having more KL
regions. In all algorithms, the inter-processor
communication cost in each CG iteration becomes a factor when the
local problems are relatively cheap. A balancing preconditioner for
the Darcy region \cite{Mandel, PenchevaYotov2003} could be used to
reduce the number of iterations. A possible extension of this work is
to incorporate stochasticity in the source term of the Stokes equation,
e.g., in modeling rainfall in surface-subsurface flow simulations. In
that case, one would expect even larger computational savings from the
use of the stochastic multiscale flux basis.

\bibliographystyle{abbrv}
\bibliography{stoch-msbasis-SD}

\end{document}